\documentclass[10pt,twoside]{article}

\usepackage[english,catalan,galician,activeacute]{babel}
\usepackage[T1]{fontenc}
\usepackage{inputenc}
\usepackage{amsmath,bm}
\usepackage{amssymb}
\usepackage{amsxtra}
\usepackage{theorem}
\usepackage{epic}
\usepackage{eepic}
\usepackage{fancybox}
\usepackage{boxedminipage}
\usepackage{calc}
\usepackage{ifthen}
\usepackage{paralist}
\usepackage{exscale,relsize}
\usepackage[rm,sc,small,center]{titlesec}
\titlelabel{\thetitle.\enspace}

\newcommand{\bdd}[2]{ {\cal B} ( #1 , #2 )}

\newcommand{\liss}[2]{ {\cal L}is ( #1 , #2 )}

\newcommand{\bip}[1]{ {\cal BIP} \left( #1 \right)}

\newcommand{\hinfty}[1]{ {\cal H}^{\infty} \left( #1 \right) }

\newcommand{\lp}[2]{ {\rm L}_{ #1 } ( #2 ) }

\newcommand{\sob}[3]{ {\rm H}^{\scriptscriptstyle #1 }_{ #2 } ( #3 ) }
\accentedsymbol{\hcirc}{ {\overset{\scriptscriptstyle \circ }{ {\rm H} }}}

\newcommand{\sobb}[3]{ {\rm W}^{\scriptscriptstyle #1 }_{ #2 } ( #3 ) }

\newcommand{\cont}[3]{ {\rm C}^{\scriptscriptstyle #1 }_{ #2 } ( #3 ) }

\newcommand{\zero}[1]{ \mbox{}_{\scriptscriptstyle 0}#1 }

\newcommand{\norm}[2]{ \| #1 \|_{ #2 } }
\newcommand{\ol}[1]{ \overline{ #1 } }

\newcommand{\eproof}{ \hfill $\Box$ }
\newcommand{\N}{\mathbb{N}}

\newcommand{\R}{\mathbb{R}}

\newcommand{\sk}[2]{ ( #1 | #2  ) }
\newcommand{\ska}[3]{ ( #1 | #2 )_{ #3 } }
\newcommand{\du}[2]{ \langle #1 , #2  \rangle }

\newcommand{\divv}{ \nabla \!\! \cdot \!}

\newcommand{\del}[2]{ \partial^{#1}_{#2} }

\newtheorem{theorem}{Theorem}[section]
\newtheorem{lemma}{Lemma}[section]

\newtheorem{proposition}{Proposition}[section]
{\theorembodyfont{\rmfamily} }

\begin{document}
\selectlanguage{english}
\thispagestyle{empty}
\numberwithin{equation}{section}
\begin{center}
{\LARGE\bf
Strong solutions in the dynamical theory of compressible fluid mixtures
}
\end{center}
\vspace{0.5cm}
\begin{center}
Matthias Kotschote
\\
{\small Fachbereich Mathematik und Statistik, Universit\"at Konstanz,
\\
D-78457 Konstanz, Germany}
\\
{\footnotesize Email: matthias.kotschote@uni-konstanz.de}
\end{center}
\vspace{0.5cm}
\begin{center}
Rico Zacher
\\
{\small Martin-Luther-Universit\"at Halle-Wittenberg,
\\
Institut f\"ur Mathematik
\\Theodor-Lieser-Stra\ss e 5, D-06120 Halle (Saale), Germany}
\\
{\footnotesize Email: rico.zacher@mathematik.uni-halle.de}
\end{center}
\vspace{0.5cm}
\begin{abstract}
In this paper we investigate the compressible Navier-Stokes-Cahn-Hilliard
equations (the so-called NSCH model) derived by Lowengrub and Truskinowsky.
This model describes the flow of a binary compressible mixture; the fluids
are supposed to be macroscopically immiscible, but partial mixing is
permitted leading to narrow transition layers.
The internal structure and macroscopic  dynamics of these layers are
induced by a Cahn-Hilliard law that the mixing ratio satisfies.
The PDE constitute a strongly coupled hyperbolic-parabolic
system. We establish a local existence and uniqueness result for
strong solutions.
\end{abstract}
\mbox{}
\section{Introduction and main result}
In this paper we are concerned with the compressible
Navier-Stokes-Cahn-Hilliard equations (the so-called NSCH model)
derived by Lowengrub and Truskinowsky, see
(\ref{eq:nsch-1})-(\ref{bc:3}) below. This system is a diffuse
interface model for the flow of a binary mixture of compressible,
viscous, and macroscopically immiscible fluids. Our main objective
is to prove existence and uniqueness of local (in time) strong
solutions of this system. Before giving the precise mathematical
formulation and stating our main result we provide some physical
background of the model.

One way to describe the flow of immiscible fluids and the motion of
interfaces between these fluids is based on the assumption that
Euler or Navier-Stokes equations apply to both sides of the
interface and across this interface certain jump conditions are
prescribed. However such a model breaks down when near interfaces a
molecular mixing of the immiscible fluids occurs in such a large
amount that the model of sharp interfaces cannot be maintained.
Another problem of such models concerns the description of merging
and reconnecting interfaces. One way out is to replace the sharp
interface by a narrow transition layer, that is one allows a
partial mixing in a small interfacial region.

For this purpose one first introduces the mass concentrations $c_i =
M_i/M$ with $M=M_1+M_2$, where $M_i$ denotes the mass of the fluid
$i$ in the representative volume $V$. Notice that this implies $c_1
+ c_2 = 1$ as well as $0 \le c_i \le 1$. A basic hypothesis is the
identification of an order parameter $c$ with a constituent
concentration, e.g. $c=c_1$, or with the difference of both
concentrations, $c = c_1 - c_2 \equiv 2 c_1 - 1$. Choosing the
latter case, $c$ varies continuously between $-1$ and $1$ in the
interfacial region and takes the values $-1$ and $1$ in the absolute
fluids. Let $u_1$, $u_2$ denote the velocities of the corresponding
fluids and $\tilde{\rho}_1 := \frac{M_1}{V}$, $\tilde{\rho}_2 :=
\frac{M_2}{V}$ the associated apparent densities which both fulfil
the equation of mass balance. Then, introducing the total density
$\rho := \tilde{\rho}_1 + \tilde{\rho}_2$ and the mass-averaged
velocity $\rho u := \tilde{\rho}_1 u_1 + \tilde{\rho}_2 u_2$, we
obtain the equation of mass balance for $\rho$ and $u$, 
\begin{align*}
\del{}{t} \rho + \divv ( \rho u) = 0, \quad (t,x) \in J \times \Omega.
\end{align*}
The total energy $E_G(t)$ in a volume $G \subset \Omega$ is to be given as the sum of
kinetic energy and (specific) Helmholtz free energy, that is it is assumed that
\begin{align*}
E_G(t) := \int\limits_G  \frac{1}{2} \rho |u|^2 + \rho \psi(\rho,c,\nabla c) \, dx.
\end{align*}
Here $\psi$ denotes the specific Helmholtz free energy density at a given temperature,
which may depend on $\rho$, $c$ and $\nabla c$. If we choose $\psi(\rho,c,\nabla c)$ as
follows
\begin{align*}
\psi(\rho,c,\nabla c) := \ol{\psi}(\rho,c) + \tfrac{1}{2} \varepsilon(\rho,c) |\nabla c|^2,
\end{align*}
also being known as the Cahn-Hilliard specific free energy density,
then the convected analogue of the Cahn-Hilliard equation can be
derived (using the second law of thermodynamics/local dissipation
inequality etc., see \cite{lowen-trusk}), that is
\begin{align*}
\del{}{t} (\rho c) + \divv( \rho u c) = \divv ( \gamma(\rho,c) \nabla \mu ),
\quad (t,x) \in J \times \Omega.
\end{align*}
The generalised chemical potential $\mu$ is given by
\begin{equation*}
\begin{split}
\mu & = \del{}{c} \psi - \rho^{-1}
\divv \left( \rho \frac{ \del{}{} \psi}{ \del{}{} \nabla c} \right)
\equiv \del{}{c} \psi - \rho^{-1} \divv \left( \rho \varepsilon(\rho,c) \nabla c \right),
\quad \del{}{c} \psi = \ol{\psi}_c(\rho,c)
+ \tfrac{1}{2} \varepsilon_c(\rho,c) |\nabla c|^2.
\end{split}
\end{equation*}
Here the parameter $\varepsilon(\rho,c)>0$ measures the interfacial
thickness and $\gamma(\rho,c) > 0$ the mobility of  the
concentration field $c$. Further, it is supposed that the stress
tensor $\mathcal{T}$ is given as the sum of a viscous and
non-viscous contribution, that is $\mathcal{T} := \mathcal{S} +
\mathcal{P}$ with
\begin{align*}
\mathcal{S}(\rho,c,u) := 2 \eta(\rho,c) \mathcal{D}(u)
+
\lambda(\rho,c) \divv u \, \mathcal{I},
\quad
\mathcal{D}(u) := \tfrac{1}{2} ( \nabla u + \nabla u^T ),
\end{align*}
where $\mathcal{I}$ denotes the identity, $\mathcal{S}$ the Cauchy stress tensor with
viscosity coefficients $\eta(\rho,c)$ and $\lambda(\rho,c)$, and
$\mathcal{P}$ the non-hydrostatic Cauchy stress tensor, which is assumed to be
of the form
\begin{align*}
\mathcal{P}(\rho,c) & := - \rho^2 \del{}{\rho} \psi \, \mathcal{I}
-
\rho \nabla c \otimes \frac{ \del{}{} \psi}{ \del{}{} \nabla c}
=
- \rho^2 \del{}{\rho} \psi \, \mathcal{I} - \rho \varepsilon(\rho,c) \nabla c \otimes \nabla c,
\\
\del{}{\rho} \psi & = \del{}{\rho} \ol{\psi} +
\tfrac{1}{2} \varepsilon_\rho(\rho,c) |\nabla c|^2.
\end{align*}
The given function $\pi := \rho^2 \psi_{\rho}$
constitutes the pressure and the extra contribution
$- \rho \nabla c \otimes \frac{ \del{}{} \psi}{ \del{}{} \nabla c}$
in the stress tensor represents capillary forces due to surface tension.
Thus the Navier-Stokes equations read as
\begin{align*}
\del{}{t} (\rho u) + \divv ( \rho u \otimes u )
-
\divv ( \mathcal{S} + \mathcal{P} ) = \rho f_{ext},
\quad (t,x) \in J \times \Omega,
\end{align*}
where $f_{ext}$ stands for external forces. A complete derivation
of this model can be found in \cite{lowen-trusk}, cf. also
\cite{gurtin} and \cite{af}.

We point out that the basic energy identity is obtained by
multiplying the momentum equation in \eqref{eq:nsch-1} by $u$,
integrating over $\Omega$, integration by parts, and using the
identity $\divv \mathcal{P} \equiv - \rho \nabla ( \psi + \rho
\del{}{\rho} \psi ) + \rho \, \mu \nabla c$. This leads to the
result
\begin{equation} \label{eq:energy}
\frac{d}{dt} E_{\Omega} (t) + \int\limits_\Omega
\mathcal{S} : \mathcal{D} \, dx + \int\limits_\Omega
\gamma(\rho,c) |\nabla \mu|^2 \, dx = \int\limits_\Omega \rho
f_{ext} \cdot u \, dx, \quad \forall t > 0.
\end{equation}

To become more specific, let $\Omega \subset \R^n$ be a bounded
domain with boundary $\Gamma := \partial \Omega$ of class $C^{4}$
decomposing as $\Gamma = \Gamma_d \cup \Gamma_s$ with dist$( \Gamma_d, \Gamma_s)>0$, where
one of these sets may be empty. The outer unit normal of $\Gamma$ at
position $x$ is denoted by $\nu(x)$. Further, let $J=[0,T]$ be a
compact time interval. The two-component (binary) viscous
compressible fluid is characterized by its total density (of the
mixture) $\rho: J \times \ol{\Omega} \to \R_+$, its mean velocity
field $u: J \times \ol{\Omega} \to \R^n$, and the mass concentration
difference of the two components (the order parameter) $c: J \times
\ol{\Omega} \to \R$. Collecting the equations from above, the
unknown functions $\rho$, $u$, and $c$ are governed by the
Navier-Stokes-Cahn-Hilliard (NSCH) system
\begin{equation} \label{eq:nsch-1}
\begin{aligned}
\del{}{t} (\rho u) + \divv ( \rho u \otimes u ) - \divv \mathcal{S}
- \divv \mathcal{P} & = \rho f_{ext},
& \quad (t,x) \in J \times \Omega,
\\
\del{}{t} ( c \rho ) + \divv ( c \rho u ) - \divv ( \gamma(\rho,c) \nabla \mu ) & = 0,
& \quad (t,x) \in J \times \Omega,
\end{aligned}
\end{equation}
\begin{equation} \label{eq:me}
\quad \del{}{t} \rho + \divv( \rho u) = 0, \quad (t,x) \in J \times \Omega,
\end{equation}
with
\begin{equation} \label{law:1}
\begin{split}
\mathcal{S} & = 2 \eta(\rho,c) \mathcal{D}(u) + \lambda(\rho,c) \divv u \, \mathcal{I},
\quad
\mathcal{P} = - \big( \pi + \tfrac{1}{2} \rho^2 \varepsilon_\rho(\rho,c) |\nabla c|^2 \big) \mathcal{I}
- \rho \varepsilon(\rho,c) \nabla c \otimes \nabla c,
\\
\mu & = \del{}{c} \psi - \rho^{-1} \nabla \cdot ( \varepsilon(\rho,c) \rho \nabla c ),
\quad
\psi  = \ol{\psi}(\rho,c) + \tfrac{1}{2} \varepsilon(\rho,c) |\nabla c|^2,
\quad
\pi = \rho^2 \del{}{\rho} \ol{\psi}.
\end{split}
\end{equation}
These equations have to be complemented by initial conditions
\begin{equation} \label{ic}
u(0,x) = u_0(x), \quad
c(0,x) = c_0(x), \quad
\rho(0,x) = \rho_0(x), \quad x \in \Omega,
\end{equation}
and boundary conditions. Two natural boundary conditions are of interest
for $u$, namely the non-slip condition
\begin{align} \label{bc:1}
u = 0, \quad (t,x) \in J \times \Gamma_d
\end{align}
and the pure slip condition
\begin{align} \label{bc:2}
\ska{ u }{ \nu }{} = 0, \quad
\mathcal{Q} \mathcal{S} \cdot \nu
=
2 \eta(\rho,c) \mathcal{Q} \mathcal{D}(u) \cdot \nu = 0,
\quad
(t,x) \in J \times \Gamma_s
\end{align}
with $\mathcal{Q}(x) := \mathcal{I} - \nu(x) \otimes \nu(x)$.
As boundary conditions for $c$, we consider
\begin{align} \label{bc:3}
\del{}{\nu} \mu(\rho,c)(t,x) = 0, \quad \del{}{\nu} c(t,x) =0,
\quad (t,x) \in J \times \Gamma,
\end{align}
meaning that no diffusion through the boundary occurs and the diffuse interface is
orthogonal to the boundary of the domain.

Note that problem \eqref{eq:nsch-1}-\eqref{bc:3} has a quasilinear
structure, since among others $\rho$ is present in front of
$\del{}{t} u$ and $\del{}{t} c$.

We are looking for strong solutions in the ${\rm L}_p$-setting. More
precisely, we seek solutions $(u,c,\rho)\in \mathcal{Z}(J):=
\mathcal{Z}_1 \times \mathcal{Z}_2  \times \mathcal{Z}_3$ where
\begin{equation} \label{spaces:Z}
\begin{split}
\mathcal{Z}_1 & := \sob{3/2 }{p}{J;\lp{p}{\Omega;\R^n}} \cap
\sob{1}{p}{J;\sob{2}{p}{\Omega;\R^n}} \cap
\lp{p}{J;\sob{4}{p}{\Omega;\R^n}},
\\
\mathcal{Z}_2  & := \sob{1}{p}{J;\lp{p}{\Omega}} \cap
\lp{p}{J;\sob{4}{p}{\Omega}}, \\
\mathcal{Z}_3 & := \sob{2 + 1/4 }{p}{J;\lp{p}{\Omega}} \cap
\cont{1}{}{J;\sob{2}{p}{\Omega}}\cap
\cont{}{}{J;\sob{3}{p}{\Omega}},
\end{split}
\end{equation}
and
\begin{align} \label{cond:ps}
p \in (\hat{p}, \infty), \quad \hat{p} := \max \left \{ 4, n \right \}.
\end{align}
Here and in the sequel the symbols ${\rm H}^s_p$ and ${\rm W}^s_p$
refer to Bessel potential spaces and Sobolev-Slobodeckij spaces,
respectively. We also write $\mathcal{Z}_i(J)$ and $\mathcal{Z}(J)$
to indicate the time interval.

To motivate the chosen solution class, let us first consider the
equation for $c$ in the base space $\lp{p}{J;\lp{p}{\Omega}}$, which
is a natural choice when looking for strong solutions. Since $\mu$
contains second order derivatives of $c$ w.r.t.\ the spatial
variables, the equation is of fourth order in space and hence
$\mathcal{Z}_2$ is the natural regularity class for $c$. Observe
that the Cahn-Hilliard equation contains a third order term of
$\rho$, that is we need $\rho\in \lp{p}{J;\sob{3}{p}{\Omega}}$ at
least. Since $\rho$ is governed by the hyperbolic equation
(\ref{eq:me}), there is no gain of regularity, that is $u\in
\lp{p}{J;\sob{4}{p}{\Omega;\R^n}}$ is required. To obtain this
regularity for the velocity, we are in turn forced to study the
Navier-Stokes equation in the base space
$\lp{p}{J;\sob{2}{p}{\Omega;\R^n}}$ at least. Note that if $c\in
\mathcal{Z}_2$ we have due to the mixed derivative theorem (cf.\
\cite{pruess-1})
\begin{align*}
c \in \mathcal{Z}_2=\sob{1}{p}{J;\lp{p}{\Omega}} \cap
\lp{p}{J;\sob{4}{p}{\Omega}}\hookrightarrow
\sob{1/2}{p}{J;\sob{2}{p}{\Omega}},
\end{align*}
and thus the natural regularity class for $\divv \mathcal{P}$, which
contains second order terms of $c$, is the space
\begin{align*}
\mathcal{X}_1^n(J) := \sob{1/2}{p}{J;\lp{p}{\Omega;\R^n}} \cap
\lp{p}{J;\sob{2}{p}{\Omega;\R^n}}.
\end{align*}
Considering $\divv \mathcal{P}$ as input for the Navier-Stokes
equation, that is taking $\mathcal{X}_1^n(J)$ as the base space for
this equation one expects that $u$ belongs to the space
$\mathcal{Z}_1$, from the maximal ${\rm L_p}$-regularity point of
view. Finally, once we have $u\in \mathcal{Z}_1$ the continuity
equation yields $\rho\in \mathcal{Z}_3$ as we will show below. Note
that the Navier-Stokes and the Cahn-Hilliard equation are strongly
coupled as $\partial_t u$, $\divv \mathcal{S}(u)$, and $\divv
\mathcal{P}$ are of the same order. The condition (\ref{cond:ps}) on
$p$ ensures the validity of several embeddings which are needed for
deriving suitable estimates for the nonlinear terms.

Our main result on the system \eqref{eq:nsch-1}-\eqref{bc:3} is the following.
\begin{theorem} \label{theo:main:1}
Let $\Omega\subset \R^n$ be a bounded domain with compact
$C^{4}$-boundary $\Gamma$ decomposing disjointly as $\Gamma =
\Gamma_d \cup \Gamma_s$ with dist$\,(\Gamma_d,\Gamma_s)>0$, $J_0=[0,T_0]$ with $T_0 \in (0,\infty)$, and $p
\in (\hat{p},\infty)$. Let further the following assumptions be satisfied.
\begin{itemize}
\item[(i)] $\varepsilon$, $\eta$, $\lambda \in \cont{4}{}{\R^2}$, $\ol{\psi}\in \cont{5}{}{\R^2}$, $\gamma \in \cont{2}{}{\R^2}$;
\item[(ii)] $\eta$, $2\eta+\lambda$, $\varepsilon$, $\gamma>0$ in $\R^2$;
\item[(iii)] $f_{ext}\in \mathcal{X}_1^n(J_0)=\sob{1/2}{p}{J_0;\lp{p}{\Omega;\R^n}} \cap
\lp{p}{J_0;\sob{2}{p}{\Omega;\R^n}}$;
\item[(iv)] $(u_0,c_0,\rho_0)\in \mathcal{V} := \sobb{4- \frac{2}{p} }{p}{\Omega;\R^n} \times \sobb{4
- \frac{4}{p} }{p}{\Omega} \times \{ \varphi \in \sob{3}{p}{\Omega}:
\varphi(x) > 0, \: \forall x \in \ol{\Omega} \}$;
\item[(v)] the subsequent compatibility conditions hold:
\begin{align*}
& u_{0|\Gamma_d} = 0, \quad \ska{ u_0 }{ \nu }{|\Gamma_s} = 0, \quad
\mathcal{Q} \mathcal{S}_{|t=0,\Gamma_s} \cdot \nu_{|\Gamma_s} =
0, \quad \del{}{\nu} c_{0} =0, \quad \del{}{\nu} \mu(\rho_0,c_0) =
0,
\\
& - \divv \mathcal{S}_{|t=0,\Gamma_d} = ( \divv \mathcal{P} + \rho
f_{ext} )_{|t=0, \Gamma_d} \in
\sobb{2-\frac{3}{p}}{p}{\Gamma_d;\R^n},
\\
& - \ska{ \divv \mathcal{S}_{|t=0} }{ \nu }{|\Gamma_s} = \ska{ \divv
\mathcal{P} - \rho \nabla u \cdot u + \rho f_{ext}}{ \nu
}{|t=0,\Gamma_s} \in \sobb{2-\frac{3}{p}}{p}{\Gamma_s},
\\
& - \mathcal{Q} \mathcal{S}( \divv \mathcal{S} )_{|t=0}  \cdot
\nu_{|\Gamma_s} = \mathcal{Q} \mathcal{S} ( \divv \mathcal{P} - \rho
\nabla u \cdot u + \rho f_{ext} )_{|t=0,\Gamma_s}  \cdot
\nu_{|\Gamma_s} \in \sobb{1-\frac{3}{p}}{p}{\Gamma_s;\R^n}.
\end{align*}
\end{itemize}
Then the system \eqref{eq:nsch-1}-\eqref{bc:3} possesses a unique
strong solution $w=(u,c,\rho)$ on a maximal time interval $J_*:=[0,T^*)$,
$T^* \le T_0$ if the solution is not global; the solution $w$ belongs to the class
$\mathcal{Z}(J_1)$ for each interval $J_1 = [0,T_1]$ with $T_1 < T^*$, or to the class
$\mathcal{Z}(J_0)$ if the solution exists globally.
The maximal time interval $J_*$ is characterized by the property:
\begin{align} \label{cond:blow-up}
\lim_{t \to T^*} w(t) \quad \text{does not exist in $\mathcal{V}_p$,}
\end{align}
where $\mathcal{V}_p$ is defined as the space of all $(u_1,c_1,\rho_1)\in \mathcal{V}$ such that the
compatibility conditions in (v) hold with $(u_0,c_0,\rho_0)$ being
replaced by $(u_1,c_1,\rho_1)$. Moreover, for fixed $f_{ext}$ not depending on $t$
the solution map $w_0 \mapsto  w(\cdot)$ generates a
local semiflow on the phase space $\mathcal{V}_{p}$.
\end{theorem}
Our result is on the original Lowengrub-Truskinovsky system. A
similar model has recently been studied by Abels and Feireisl
\cite{af}. They proved existence of global weak solutions, but not
uniqueness, for a simplified version of the Lowengrub-Truskinovsky
system where the Helmholtz free energy (in our notation) is given by
\begin{align*}
F = \int\limits_\Omega  \big( \rho
\ol{\psi}(c,\rho)+\frac{1}{2}\,|\nabla c|^2 \big) \, dx,
\end{align*}
that is $\varepsilon=1/\rho$, see also Anderson \cite[p. 151]{AMFW}. The
approach in \cite{af} does not seem to extend to the original
Lowengrub-Truskinovsky system, since the energy estimates, on which
the method in \cite{af} is based, do not provide any bound for
$\nabla c$ in vacuum zones, i.e. where $\rho=0$. A similar model for
incompressible fluids was studied by Boyer \cite{boyer}, Liu and
Shen \cite{liu-shen}, Starovoitov \cite{staro}, and Abels
\cite{abels-1}.

The basic tool in our proof of Theorem \ref{theo:main:1} is the
contraction mapping principle. We proceed as follows. Regarding
$u$ as given and assuming sufficient regularity, the continuity
equation, as well-known in the literature, can be solved by means
of the method of characteristics provided that the condition
$\ska{ u }{ \nu }{|\Gamma } \ge 0$ is satisfied, see e.g.\
\cite{s1}. Inserting the solution $\rho=L[u] \rho_0$ into the
equations for $u$ and $c$ reduces the original system to a
non-local, fully nonlinear, strongly coupled system for $u$ and
$c$. This problem then is locally solved by means of a fixed point
argument using maximal regularity for the linearized problem. The
unique solution $(u,c)$ is found in the class $\mathcal{Z}_1\times
\mathcal{Z}_2$, and this in turn gives rise to $\rho\in
\mathcal{Z}_3$ via the relation $\rho=L[u] \rho_0$.

When looking at the reduced problem for $(u,c)$ one realizes that
it is impossible to derive a contraction inequality in the space
$\mathcal{Z}_1\times \mathcal{Z}_2$. To overcome this difficulty
we will work with a larger base space of the fixed point mapping.
This idea was independently introduced by Kato and Lax in the context
of quasilinear hyperbolic symmetric systems and has already often been used
in the literature, see e.g. \cite{M1}, \cite{s1}, \cite{gkz}. It turns out that the space
$Z_1(J)\times Z_2(J)$ defined by
\begin{align*}
Z_1(J) & := \sob{5/4}{2}{J;\lp{2}{\Omega;\R^n}} \cap
\sob{1}{2}{J;\sob{1}{2}{\Omega;\R^n}} \cap \lp{2}{J;\sob{3}{2}{\Omega;\R^n}},
\\
Z_2(J) & := \sob{3/4}{2}{J;\lp{2}{\Omega}} \cap
\lp{2}{J;\sob{3}{2}{\Omega}}
\end{align*}
is a suitable choice for the contraction property.

There is still another problem that arises in setting up the fixed point argument,
it is due to the nonlinear boundary condition $\partial_\nu \mu(\rho,c)=0$. It seems that for the derivation of the
desired contraction property one requires that the identity $\partial_\nu \mu=0$ is preserved under the fixed
point mapping. To overcome this difficulty we add the variable $\mu$, that is we work with
triples $(u,c,\mu)$ and view the reduced problem for $(u,c)$ as a problem for $(u,c,\mu)$.
Since $c \in \mathcal{Z}_2(J)$ and
\begin{align*}
\mu=\del{}{c} \psi - \rho^{-1} \divv \left( \varepsilon \rho \nabla c \right),
\end{align*}
the natural regularity class of $\mu$ is given by
\begin{align*}
\mu \in \mathcal{Z}_\mu(J) := \sob{1/2}{p}{J;\lp{p}{\Omega}} \cap
\lp{p}{J;\sob{2}{p}{\Omega}}.
\end{align*}
As differences of $c$ are considered in $Z_2(J)$, the natural space for the contraction estimate for $\mu$ will be
\begin{align*}
Z_\mu(J) :=\sob{1/4}{2}{J;\lp{2}{\Omega}} \cap \lp{2}{J;\sob{1}{2}{\Omega}}.
\end{align*}
We put
\begin{equation} \label{mu0def}
\mu_0:=\mu_{|t=0}= \ol{\psi}_c(c_0,\rho_0)
+ \tfrac{1}{2} \varepsilon_c(c_0,\rho_0) |\nabla c_0|^2  - \rho_0^{-1} \divv \left( \varepsilon(c_0,\rho_0) \rho_0 \nabla c_0 \right).
\end{equation}
Let further $u_\bullet$ be the trace of $\partial_t u$ at $t=0$ which is obtained from \eqref{eq:nsch-1}-\eqref{ic}.  Note that
\begin{equation*}
\partial_t u\in \sob{1/2}{p}{J;\lp{p}{\Omega;\R^n}} \cap \lp{p}{J;\sob{2}{p}{\Omega;\R^n}} 
\hookrightarrow \cont{}{}{J;\sobb{2-\frac{4}{p}}{p}{\Omega;\R^n}},
\end{equation*}
see e.g.\ \cite{zacher}, so that $\sobb{2-\frac{4}{p}}{p}{\Omega;\R^n}$ is the natural space for $u_\bullet$. 
For $T \in (0,T_0]$ (with $T_0>0$ being fixed) we set $J:=[0,T]$ and will consider the set
\begin{equation} \label{set:SIGMA}
\begin{split}
\Sigma_{T} := \{ &(u,c,\mu) \in \mathcal{Z}_1(J) \times
\mathcal{Z}_2(J)\times \mathcal{Z}_\mu(J):\;\,(u,\partial_t u,c,\mu)(0) = (u_0,u_\bullet,c_0,\mu_0),\\
& \;u=0\;\mbox{on}\;J\times \Gamma_d,\;\ska{ u }{ \nu }{} = 0\;\mbox{and}\;
\mathcal{Q} \mathcal{D}(u) \cdot \nu=0\;\mbox{on}\;J\times\Gamma_s,
\\
& \;\partial_\nu c=\partial_\nu \mu=0\;\mbox{on}\;J\times \Gamma,\;\mbox{and}\\
& \;
\norm{ (u,c,\mu) - (\ol{u},\ol{c},\ol{\mu}) }{
{\mathcal{Z}_1}(J) \times {\mathcal{ Z}_2}(J)\times \mathcal{Z}_\mu(J) } \le 1 \},
\end{split}
\end{equation}
in which local solutions $(u,c,\mu)$ of the reduced problem will be
sought. Here $(\ol{u},\ol{c},\ol{\mu})\in \mathcal{Z}_1(J_0) \times
\mathcal{Z}_2(J_0)\times\mathcal{Z}_\mu(J_0) $, with $J_0=[0,T_0]$, is a triple of certain
reference functions with $(\ol{u},\partial_t \ol{u},\ol{c},\ol{\mu})(0)=(u_0,u_\bullet,c_0,\mu_0)$. We will
show that for sufficiently small $T$ the fixed point mapping
associated with the reduced problem leaves $\Sigma_T$
invariant and is a strict contraction in $Z_1(J)\times Z_2(J)\times Z_\mu(J)$.
Therefore the contraction mapping principle applies, since
$\Sigma_T$ is a closed subset of $Z_1(J)\times Z_2(J)\times Z_\mu(J)$, see Lemma \ref{lem:1} below.

The paper is organized as follows. In Section 2 we fix some
notation and provide some auxiliary results. In Section
\ref{sec:3} we describe our reformulation of the system
\eqref{eq:nsch-1}-\eqref{bc:3} as a fixed point problem. Section
\ref{sec:MR} is devoted to the maximal regularity property for the
linearized problem. In particular we establish maximal regularity
in the non-standard higher regularity class $\mathcal{Z}_1(J)$ for
the subproblem for $u$, i.e. the linearized Navier-Stokes problem.
In Section \ref{sec:CE} we study the continuity equation. We prove
regularity results for $\rho$ given a velocity
$u\in\mathcal{Z}_1(J)$ and derive a contraction estimate in the space
\begin{align*}
Z_3(J) & := \sob{2}{2}{J;\lp{2}{\Omega}} \cap \cont{1}{}{J;\sob{1}{2}{\Omega}} \cap \cont{}{}{J;\sob{2}{2}{\Omega}}.
\end{align*}
In Section \ref{sec:WCH} we prove the crucial contraction estimate
for the Cahn-Hilliard subproblem in the space $Z_2(J)\times Z_\mu(J)$. The proof
of Theorem \ref{theo:main:1} is completed in Section
\ref{sec:PrThm}, where we make use of the estimates from the
preceding sections to carry out the fixed point argument. The
paper concludes with generalizations of the main result (Section
\ref{sec:8}).

\section{Preliminaries} \label{sec:2}
We begin by fixing some notation. If $X$ is a Banach space and $\Omega$ is a Lebesgue
measurable subset of $\R^n$, then for $1< p< \infty$ and $s> 0$ the symbols ${\rm H}^s_p(\Omega;X)$ 
and ${\rm W}^s_p(\Omega;X)$ stand for the Bessel potential spaces resp.\ Sobolev-Slobodeckij spaces 
of $X$-valued functions on $\Omega$. For the Lebesgue spaces and spaces of continuous or H\"older 
continuous functions we use standard notation. Furthermore, if
$\mathcal{F}(I;X)$ is any function space with $I=[0,T] \subseteq
\R_+$, we set $\zero{ \mathcal{F}(I;X) } := \{ v \in
\mathcal{F}(I;X): (\partial^k_t v)_{|t=0} = 0,\,k=0,1,2,\ldots,\;\mbox{whenever the trace exists}\}$ and  
${}^0{ \mathcal{F}(I;X) } := \{ v \in \mathcal{F}(I;X): (\partial^k_t v)_{|t=T} = 0,\,k=0,1,2,\ldots,\;\mbox{whenever the trace exists}\}$. 
Recall that a Banach space $X$ belongs to the class ${\cal HT}$, if the Hilbert transform is bounded on $L_2(\R;X)$.

Let $X$ be a complex Banach space and $A$ be a closed linear
operator in $X$. Then $A$ is called {\em pseudo-sectorial} if
$(-\infty,0)$ is contained in the resolvent set of $A$ and the
resolvent estimate $|t(t+A)^{-1}|_{{\cal B}(X)}\le C,\, t>0$,
holds, for some constant $C>0$. If in addition the null space
${\cal N}(A)=\{0\}$, and the domain ${\cal D}(A)$ as well as the
range ${\cal R}(A)$ of $A$ are dense in $X$ then $A$ is called
{\em sectorial}. The class of sectorial operators in $X$ is
denoted by ${\cal S}(X)$, and by $\phi_A$ we mean the spectral angle of $A\in {\cal S}(X)$. For $1\le p< \infty$, 
and $\gamma\in(0,1)$ the space $(X,D_A)_{\gamma,\,p}$ denotes the real interpolation
space $(X,D_A)_{\gamma,\,p}$, where $D_A$ stands for the domain of
$A$ equipped with the graph norm. We say that an operator $A\in {\cal S}(X)$ admits bounded
imaginary powers and write $A\in {\cal BIP}(X)$ if the imaginary
powers $A^{is}$ form a bounded $C_0$-group on $X$. The type
$\theta_A$ of this group is called the {\em power angle} of $A$; there holds $\theta_A\ge \phi_A$. 
We refer to \cite{dhp} for properties of operators from
the classes  ${\cal S}(X)$ and ${\cal BIP}(X)$ (and other important subclasses of ${\cal S}(X)$).

The following proposition can be found in \cite{pruess-1}.
\begin{proposition} \label{prop:1}
Let $1 < p < \infty$, $1/p < \beta < 1$, suppose $A$ is an invertible pseudo-sectorial
operator in $X$ with $\phi_A < \pi/2$, and set $u(t) = e^{-A t} x$, $x \in X$. Then the
following statements are equivalent:
\begin{align*}
(i) \quad x \in D_A(\beta - 1/p,p); \quad (ii) \quad u \in \lp{p}{\R_+;D_A(\beta,p)}; \quad
(iii) \quad u \in \sobb{\beta}{p}{\R_+;X}.
\end{align*}
\end{proposition}

Next we collect several embeddings for the function spaces \eqref{spaces:Z}
which will be used frequently in what follows. By the mixed derivative theorem, see
\cite{pruess-1}, we have the embeddings
\begin{equation} \label{embedd:1}
\begin{split}
\mathcal{Z}_1 & \hookrightarrow \sob{1 +
\theta/2}{p}{J;\sob{(1-\theta)2}{p}{\Omega;\R^n}} \cap
\sob{\theta}{p}{J;\sob{(1-\theta)2 + 2}{p}{\Omega;\R^n}},
\\
\mathcal{Z}_2 & \hookrightarrow
\sob{\theta}{p}{J;\sob{4(1-\theta)}{p}{\Omega}},
\quad
\mathcal{Z}_\mu \hookrightarrow \sob{\theta/2}{p}{J;\sob{2(1-\theta)}{p}{\Omega}},
\\
\mathcal{Z}_3 & \hookrightarrow
\sob{(2+1/4)\theta}{p}{J;\sob{3(1-\theta)}{p}{\Omega}},
\quad
\theta \in (0,1).
\end{split}
\end{equation}
For $p>\hat{p}$ and $0< \beta < 1/4 - 1/p$ Sobolev embeddings imply
\begin{equation} \label{embedd:2}
\begin{split}
\mathcal{ Z}_1 &\hookrightarrow \mathcal{W}_1 :=
\cont{5/4 + \beta}{}{J;\cont{}{}{\ol{\Omega};\R^n}}\cap
\cont{3/4 + \beta}{}{J;\cont{1}{}{\ol{\Omega};\R^n}} \cap
\cont{1/4 + \beta}{}{J;\cont{2}{}{\ol{\Omega};\R^n}},
\\
\mathcal{ Z}_2 & \hookrightarrow \mathcal{W}_2 :=
\cont{1/2 + \beta}{}{J;\cont{}{}{\ol{\Omega}}} \cap
\cont{1/4 + \beta }{}{J;\cont{1}{}{\ol{\Omega}}} \cap
\cont{\beta}{}{J;\cont{2}{}{\ol{\Omega}}},
\\
\mathcal{ Z}_3 & \hookrightarrow \mathcal{W}_3 :=
\cont{1}{}{J;\cont{}{}{\ol{\Omega}}} \cap
\cont{1/2}{}{J;\cont{1}{}{\ol{\Omega}}} \cap
\cont{}{}{J;\cont{2}{}{\ol{\Omega}}},
\\
\mathcal{ Z}_\mu & \hookrightarrow \mathcal{W}_\mu :=
\cont{1/4 + \beta}{}{J;\lp{p}{\Omega}} \cap
\cont{\beta}{}{J;\cont{}{}{\ol{\Omega}}} \cap
\lp{p}{J;\cont{1}{}{\ol{\Omega}}}.
\end{split}
\end{equation}
\begin{lemma} \label{lem:1}
Under the above assumptions on $\Omega$ and $p$, the set
$\Sigma_T$ as defined as in (\ref{set:SIGMA}) is closed in the space $Z_1(J)\times Z_2(J)\times Z_\mu(J)$.
\end{lemma}
{\it Proof.} The assertion of this lemma follows from the above Sobolev embeddings ($p>\hat{p} > 2$),
and the subsequent abstract lemma.
\eproof
\begin{lemma} \label{auxlemma}
Let $X$, $Y$ be Banach spaces
with $Y \hookrightarrow X$ densely and $Y$ being reflexive. Then for any $r>0$ the ball
$B_r(0) := \{ y \in Y: \|y\|_Y \le r \}$ is closed with respect to the
topology of $X$.
\end{lemma}
A proof of Lemma \ref{auxlemma} can be found, e.g., in \cite{k1}.
%
%
%
%
%
%

\begin{lemma} \label{lem:prod:2}
Let $J=[0,T]$ be a compact interval with $T\le T_0$ and $T_0$ being fixed, $1<q<\infty$, 
$\alpha\in (\frac{1}{q},1)$, $0<\varepsilon<1-\alpha$, and $X$ be a Banach space of class $\mathcal{HT}$. 
Then for any $v\in {}_0\sob{\alpha}{q}{J;X}$ and any $\varphi\in \cont{\alpha+\varepsilon}{}{J;\mathcal{B}(X)}$
the product $\varphi v$ belongs to the space ${}_0\sob{\alpha}{q}{J;X}$ as well and there holds
\begin{equation} \label{pointmultH}
\norm{\varphi v}{{}_0\sob{\alpha}{q}{J;X}}\le C\big(\norm{\varphi}{
\cont{}{}{J;\mathcal{B}(X)}}\norm{v}{{}_0\sob{\alpha}{q}{J;X}}+
\norm{\varphi}{\cont{\alpha+\varepsilon}{}{J;\mathcal{B}(X)}}
\norm{v}{\lp{q}{J;X}}\big),
\end{equation}
where the constant $C$ is independent of $T,\varphi$, and $v$.
%
\end{lemma}
{\it Proof.} Set $g_{1-\alpha}(t)=t^{-\alpha}/\Gamma(\alpha)$, $t>0$. The Riemann-Liouville fractional derivative 
of $v$ is defined as
\begin{equation*}
\partial_t^\alpha v:=\partial_t \int_0^t g_{1-\alpha}(t-\tau)v(\tau)\,d\tau,\quad t\in J.
\end{equation*}
It has been shown in \cite{zacher}, see also \cite{zacherD}, that $v\mapsto \norm{\partial_t^\alpha v}{\lp{q}{J;X}}$ 
defines an equivalent norm for the space ${}_0\sob{\alpha}{q}{J;X}$, whenever $X$ belongs to the class $\mathcal{HT}$.
For $v$ and $\varphi$ as in the statement of the lemma we further have the product rule
\begin{equation*}
\partial_t^\alpha(\varphi v)(t)=\varphi(t)\partial_t^\alpha v(t)+\int_0^t [-g_{1-\alpha}'(t-\tau)]\big(\varphi(t)-\varphi(\tau)\big)v(\tau)\,d\tau,\quad t\in J,
\end{equation*}
which can be obtained by simple algebra assuming that $\varphi\in \cont{1}{}{J;\mathcal{B}(X)}$ and an approximation argument. A version of this
product rule can already be found in \cite[Lemma 2.2]{zacher1}. We may now estimate as follows.
\begin{align*}
\norm{\varphi v}{{}_0\sob{\alpha}{q}{J;X}} & \le C_1 \norm{\partial_t^\alpha (\varphi v)}{\lp{q}{J;X}}
\le C_1\Big(\norm{\varphi \partial_t^\alpha v}{\lp{q}{J;X}}\\
& \quad +\frac{\alpha}{\Gamma(1-\alpha)}\,\norm{\int_0^t (t-\tau)^{-\alpha-1}\big(\varphi(t)-\varphi(\tau)\big)v(\tau)\,d\tau}
{\lp{q}{J;X}}\Big)\\
& \le C_1\Big(\norm{\varphi}{\cont{}{}{J;\mathcal{B}(X)}} \norm{\partial_t^\alpha v}{\lp{q}{J;X}}\\
& \quad +\frac{\alpha}{\Gamma(1-\alpha)}
\norm{\varphi}{\cont{\alpha+\varepsilon}{}{J;\mathcal{B}(X)}}
\norm{\int_0^t (t-\tau)^{\varepsilon-1}\norm{v(\tau)}{X}\,d\tau}{\lp{q}{J}}\Big)\\
&\le C_2\Big(\norm{\varphi}{\cont{}{}{J;\mathcal{B}(X)}}
\norm{v}{{}_0\sob{\alpha}{q}{J;X}}+\norm{g_\varepsilon}{\lp{1}{J}}
\norm{\varphi}{\cont{\alpha+\varepsilon}{}{J;\mathcal{B}(X)}}
\norm{v}{\lp{q}{J;X}}\Big),
\end{align*}
which implies the desired estimate with a constant that can be chosen independent of $T\in (0,T_0]$, since $v$ has vanishing trace at $t=0$.
\eproof
\section{Formulation of the fixed point problem} \label{sec:3}
In this section we describe how the original problem is reformulated as a fixed point problem.
The basic idea is to rewrite the Navier-Stokes equations as well as the
Cahn-Hilliard equation in \eqref{eq:nsch-1} such that the left-hand
side becomes linear and the nonlinearities on the right-hand side can be estimated appropriately
to make the fixed point argument work. As already mentioned in the introduction we view the Cahn-Hilliard equation,
together with the law for $\mu$, as a system for the pair $(c,\mu)$.
We further point out that the linearisation is carried out
in such a way that the elliptic operator appearing in the Cahn-Hilliard equation
maintains its divergence structure
and can be viewed as the square of an elliptic operator. This feature
will be essential to establish the desired contraction inequality for the
Cahn-Hilliard problem in the space $Z_2(J)\times Z_\mu(J)$, see Section \ref{sec:WCH}.

The governing equations for $u$, $c$, and $\mu$ can be rewritten as
\begin{equation} \label{eq:nsch-2}
\begin{aligned}
\rho_0 \del{}{t} u + \mathcal{A}(D) u  + \mathcal{B}(D) c + \mathcal{C} \mu & = F_1(w,\rho),
& \quad & (t,x) \in J \times \Omega,
\\
\tfrac{\varepsilon_0 \rho_0}{\gamma_0} \del{}{t} c
-
\divv( \varepsilon_0 \nabla\mu ) & = F_2(w,\rho),
& \quad & (t,x) \in J \times \Omega,
\\
-\mu-\divv(\varepsilon_0 \nabla c) & = F_\mu(w,\rho), & \quad & (t,x) \in J \times \Omega,\\
u & = 0, & \quad & (t,x) \in J \times \Gamma_d,
\\
( \ska{ u }{\nu}{|\Gamma_s}, \mathcal{Q} \mathcal{D}(u) \cdot \nu_{|\Gamma_s} )
& = 0,
& \quad & (t,x) \in J \times \Gamma_s,
\\
\del{}{\nu} c = \del{}{\nu} \mu & =0,
& \quad & (t,x) \in J \times \Gamma,
\\
u = u_0, \quad c & = c_0,
& \quad & (t,x) \in \{0\} \times \Omega,
\end{aligned}
\end{equation}
where now $w=(u,c,\mu)$, $a_0 := a_{|t=0}$ for $a \in \{ \eta, \lambda, \gamma, \varepsilon, \pi \}$.
Further $\mathcal{A}(D)$ and $\mathcal{B}(D)$ are second order operators, and $\mathcal{C}$ is a
multiplication operator; they are defined by
\begin{align*}
\mathcal{A}(D) u & :=
- \divv \big( 2 \eta_0 \mathcal{D}(u) +\lambda_0 \divv u \, \mathcal{I}  \big),
\\
\mathcal{B}(D) c & :=
[ \rho_0^2  \del{}{\rho} \varepsilon_{0}  + \rho_0 \varepsilon_0] \nabla c_0 \cdot \nabla^2 c,
\\
\mathcal{C} \mu & := - \rho_0 \nabla c_0 \mu.
\end{align*}
The nonlinearities $F_1$, $F_2$, and $F_\mu$ are given by
\begin{equation} \label{def:F-data}
\begin{split}
F_1(w,\rho) & := B_1(w,\rho)u + B_2(w,\rho) c + B_3(w,\rho) \rho + B_\mu(w) \mu + B_{low}(w,\rho),
\\
F_2(w,\rho) & := \tfrac{\varepsilon_0}{\gamma_0} \big \{
\del{}{t} \big( [\rho_0 - \rho] c \big)
-
\divv( c \rho u ) - \divv \big( [\gamma_0 - \gamma] \nabla \mu \big) \big \}
+
\tfrac{\varepsilon_0^2}{\gamma_0} \nabla
(
\tfrac{\gamma_0}{\varepsilon_0}
)
\cdot \nabla\mu,
\\
F_\mu(w,\rho)  &:= \divv( [\varepsilon-\varepsilon_0] \nabla c )
+ \rho^{-1} \varepsilon \nabla \rho \cdot \nabla c - \del{}{c} \psi,
\end{split}
\end{equation}
where
\begin{equation} \label{def:B}
\begin{split}
B_{low}(w,\rho) & := - \rho \nabla u \cdot u - \del{}{c} ( \rho \psi + \rho^2 \del{}{\rho} \psi ) \nabla c
+ \rho f_{ext},
\\
B_1(w,\rho) \phi  & := (\rho_0 - \rho ) \del{}{t} \phi
-
\divv \big( 2 [\eta_0 - \eta(\rho,c) ] \mathcal{D}(\phi)
+
[\lambda_0 - \lambda(\rho,c) ] \divv \phi \, \mathcal{I} \big),
\\
B_2(w,\rho) \phi & := \big ( [ \rho_0^2 \del{}{\rho} \varepsilon_0 + \rho_0 \varepsilon_0 \nabla c_0] \nabla c_0
-  [\rho^2 \del{}{\rho} \varepsilon + \rho_0 \varepsilon_0 \nabla c_0] \nabla c \big) \nabla^2 \phi,
\\
B_3(w,\rho) \phi & := - \del{}{\rho} (  \rho^2 \del{}{\rho} \psi ) \nabla \phi,
\quad
B_\mu(w,\rho) \phi := - ( \rho_0 \nabla c_0 - \rho \nabla c ) \phi.
\end{split}
\end{equation}
Observe that $B_i(w_0,\rho_0) \equiv 0$ holds for $i=1,2,\mu$ and
the term $\del{}{c} ( \rho \psi + \rho^2 \del{}{\rho} \psi ) \nabla c$ is of lower order. Moreover,
$\del{}{\rho} (  \rho^2 \del{}{\rho} \psi )$ depends only on $\rho$, $c$, and $\nabla c$.
For the 'freezing of the coefficients' in the quasilinear terms in the Navier-Stokes equation we may use
the fixed functions $c_0$ and $\rho_0$, as by imposing $p > 4$, cf. the condition \eqref{cond:ps},
the embeddings $c_0 \in \sobb{4-4/p}{p}{\Omega} \hookrightarrow \sob{3}{p}{\Omega}
\hookrightarrow \cont{2}{}{\ol{\Omega}}$ hold.
These are necessary to ensure $\mathcal{A}(D) u \in \mathcal{X}_1^n(J)$ for given $u \in \mathcal{Z}_1(J)$.

Problem \eqref{eq:nsch-2} can be viewed as an abstract equation for $w=(u,c,\mu)$ of the form
\begin{align} \label{eq:fpe-1}
\mathcal{ L}w = (\mathcal{F}_1(w,\rho),u_0,\mathcal{F}_2(w,\rho),c_0,\mathcal{F}_\mu(w,\rho))
=: \mathcal{F}(w,\rho), \quad \rho = L[u]\rho_0,
\end{align}
where $\mathcal{ L}$ stands for the linear operator on the left-hand side of \eqref{eq:nsch-2}
with the components being ordered in such a way that we first
have the equations for $u$, then for $c$, and finally for $\mu$. Further, $L[u]\rho_0$ denotes
the solution operator to the equation of conservation of mass,
see Section \ref{sec:CE}, and $\mathcal{F}_i$, $i=1,2,\mu$, comprises the non\-linearity $F_i$ as well as
the corresponding zero boundary data,
\begin{equation} \label{def:F}
\mathcal{F}_1(w,\rho) := (F_1(w,\rho), 0, 0 ),
\quad
\mathcal{F}_2(w,\rho) := (F_2(w,\rho), 0),
\quad
\mathcal{F}_\mu(w,\rho) := (F_\mu(w,\rho), 0).
\end{equation}
We will show that equation \eqref{eq:fpe-1} defines a nonlinear mapping
$\mathcal{G}: \Sigma_T \to \mathcal{Z}_1(J) \times \mathcal{Z}_2(J)\times \mathcal{Z}_\mu(J)$
(see (\ref{set:SIGMA}) for the definition of $\Sigma_T$) according to
\begin{equation} \label{eq:fpe-mapping}
\mathcal{G}: \: w \mapsto w':= (u',c',\mu'), \quad
\mathcal{L}w' = \mathcal{F}(w,\rho), \quad w = (u,c,\mu), \quad \rho = L[u]\rho_0,
\end{equation}
which for a sufficiently small $T$ leaves $\Sigma_T$ invariant and becomes a strict contraction
with respect to the weaker topology of $Z_1\times Z_2\times Z_\mu$.
\section{Maximal regularity for Cahn-Hilliard and a viscous fluid} \label{sec:MR}
To make our fixed point argument work, we have to show that the linearized problem with the operator
$\mathcal{L}$ on the left-hand side, that is the problem
\begin{equation} \label{eq:linear}
\begin{aligned}
\rho_0 \del{}{t} u + \mathcal{A}(D)u + \mathcal{B}(D) c + \mathcal{C} \mu & = f(t,x),
& \quad & (t,x) \in J \times \Omega,
\\
\tfrac{\varepsilon_0 \rho_0}{\gamma_0} \del{}{t} c
-
\divv( \varepsilon_0 \nabla\mu )
 & = g(t,x),
& \quad & (t,x) \in J \times \Omega,
\\
-\mu-\divv(\varepsilon_0 \nabla c) & = g_\mu(t,x), & \quad & (t,x) \in J \times \Omega,\\
u & = h_d(t,x),
& \quad & (t,x) \in J \times \Gamma_d,
\\
( \sk{u}{\nu}, \mathcal{Q} \mathcal{D}(u) \cdot \nu ) &  = h_s(t,x),
& \quad & (t,x) \in J \times \Gamma_s,
\\
\del{}{\nu} c = h_1(t,x) , \quad \del{}{\nu}\mu & = h_\mu(t,x),
& \quad & (t,x) \in J \times \Gamma,
\\
u = u_0(x), \quad c & = c_0(x),
& \quad & (t,x) \in \{0\} \times \Omega,
\end{aligned}
\end{equation}
has the property of maximal regularity
in the described setting. This means we have to prove that for any right-hand side data in the natural
regularity classes and subject to the necessary compatibility
conditions the linear problem (\ref{eq:linear}) possesses a unique solution
$(u,c,\mu)\in \mathcal{Z}_1(J) \times \mathcal{Z}_2(J)\times \mathcal{Z}_\mu(J)$.
Since in the linearized system the problems for $u$ and $(c,\mu)$ decouple,
one can first solve the Cahn-Hilliard problem in \eqref{eq:linear}
and insert its solution into the terms $\mathcal{B}(D) c$ and $\mathcal{C} \mu$ which become
data for the equation of $u$, which is solved in the second step.
Thus the linearized problem reduces to the study of two separated problems.

We begin with the linear Cahn-Hilliard problem, that is with the system for $(c,\mu)$.
The corresponding maximal regularity result is the following.
\begin{theorem} \label{theo:ch}
Let $\Omega$ be a bounded domain in $\R^n$ with $C^{4}$-boundary
$\Gamma$, let $J=[0,T]$ be a compact time interval, and
$p > \max \{ 1, \frac{n}{3} \}$ with $p \not= \frac{5}{3}$, $5$ .
Further, assume that $\rho_0$, $\gamma_0$, $\varepsilon_0 \in \sob{3}{p}{\Omega}$
and $\rho_0(x)$, $\gamma_0(x)$, $\varepsilon_0(x) > 0$ for all $x \in \ol{\Omega}$.
Then the linear Cahn-Hilliard subproblem in \eqref{eq:linear} possesses a unique
solution $(c,\mu) \in \mathcal{Z}_2(J)\times \mathcal{Z}_\mu(J)$ if and only if the data
$g$, $g_\mu$, $h$, $h_\mu$, $c_0$ satisfy the following conditions
\begin{compactenum}
\item $g \in \mathcal{X}_2(J) := \lp{p}{J;\lp{p}{\Omega}}$,
$g_\mu\in \mathcal{X}_\mu(J):=\sob{1/2}{p}{J;\lp{p}{\Omega}} \cap
\lp{p}{J;\sob{2}{p}{\Omega}}$;
\item $(h_1,h_\mu) \in \mathcal{Y}_1(J) \times \mathcal{Y}_\mu(J)$,
with $\mathcal{Y}_k(J):=\sobb{1 - \frac{k}{4}-\frac{1}{4p}}{p}{J;\lp{p}{\Gamma}}
\cap \lp{p}{J;\sobb{4-k-\frac{1}{p}}{p}{\Gamma}}$, $k=1,3$,
and $\mathcal{Y}_\mu(J):=\mathcal{Y}_3(J)$;
\item $c_0 \in \sobb{4- \frac{4}{p}}{p}{\Omega}$;
\item $\del{}{\nu} c_0 = h_{1|t=0}$ in
$\sobb{3-\frac{5}{p}}{p}{\Gamma}$ for $p > \frac{5}{3}$, and
$\del{}{\nu}\mu_0 = h_{\mu|t=0}$ in
$\sobb{1-\frac{5}{p}}{p}{\Gamma}$ for $p>5$, where $\mu_0:=g_{\mu|t=0}-\nabla\cdot (\varepsilon_0 \nabla c_0)$.
\end{compactenum}
Moreover, the terms $\mathcal{B}(D)c$ and $\mathcal{C} \mu$ in the first equation of \eqref{eq:linear}
belong to $\mathcal{X}_1^n(J)$ provided that $p > \hat{p}$.
\end{theorem}

{\it Proof.} This result is well-known and follows, e.g., from \cite{dhp-2},
see also \cite{prz} and \cite{pw-1}. The last assertion is a consequence of the mixed
derivative theorem and multiplier theorems for fractional Sobolev spaces. \eproof

${}$

\noindent The remaining equations of \eqref{eq:linear} represent a linear problem for $u$,
as the terms $\mathcal{B}(D)c$ and $\mathcal{C}\mu$ are determined by the previous theorem.
The following theorem gives necessary and sufficient conditions for the unique
solvability of the $u$-problem in the maximal $L_p$-regularity class and is well-known, it follows, e.g., from the results in \cite{dhp-2}.

\begin{theorem} \label{theo:stokes:1}
Let $\Omega$ be a bounded domain in $\R^n$, $n \ge 1$, with $C^{2}$-boundary
$\Gamma$ decomposing disjointly as $\Gamma = \Gamma_d \cup \Gamma_s$
with dist$\,(\Gamma_d,\Gamma_s)>0$, $J=[0,T]$, and
$p \in (1,\infty)$ with $p \not = 3/2$, $3$.
Further, assume that $\rho_0 \in \cont{}{}{\ol{\Omega}}$, $\eta_0$, $\lambda_0 \in \cont{1}{}{\ol{\Omega}}$
and $\rho_0(x)$, $\eta_0(x)$, $2 \eta_0(x) + \lambda_0(x) > 0$ for all $x \in \ol{\Omega}$.
Then the subproblem for $u$ in \eqref{eq:linear} possesses a unique solution
\begin{align*}
u \in \sob{1}{p}{J;\lp{p}{\Omega;\R^n}} \cap \lp{p}{J;\sob{2}{p}{\Omega;\R^n}},
\end{align*}
if and only if the data  $\tilde{f} = f - \mathcal{B}(D)c - \mathcal{C} \mu$, $h_d$, $h_s$, $u_0$
satisfy the following conditions.
\begin{compactenum}
\item $\tilde{f} \in \lp{p}{J;\lp{p}{\Omega;\R^n}}$;
\item $(h_d,h_s) \in Y_{0,d}(J;\R^n) \times
Y_{0,s}(J) \times Y_{1,s}(J;\R^n)$ with $h_s := (h_{s1},h_{s2})$ and \\
$Y_{i,k}(J;E):= \sobb{(2 - i - \frac{1}{p})\frac{1}{2}}{p}{J;\lp{p}{\Gamma_k;E}}
\cap \lp{q}{J;\sobb{2-i-\frac{1}{p}}{p}{\Gamma_k;E}}$, $i=0,1$, $k=d,s$;
\item $u_0 \in \sobb{2-\frac{2}{p}}{p}{\Omega;\R^n}$;
\item $u_{0|\Gamma_d} = h_{d|t=0}$ in $\sobb{2-\frac{3}{p}}{p}{\Gamma_d;\R^n}$ if $p>3/2$;
\item $\ska{ u_{0} }{ \nu}{|\Gamma_s} = h_{s1|t=0}$ in
$\sobb{2-\frac{3}{p}}{p}{\Gamma_s}$,
$\mathcal{Q} \mathcal{D}(u_0) \cdot \nu_{|\Gamma_s} = h_{s2|t=0}$ in
$\sobb{1-\frac{3}{p}}{p}{\Gamma_s;\R^n}$ if $p > 3$.
\end{compactenum}
\end{theorem}
The (non-standard) higher maximal regularity results for $u$ we need read as follows.
Here we use the notation
\begin{align*}
\mathcal{B}_d u=u_{|\Gamma_d},\quad
\mathcal{B}_s u=\big( \sk{u}{\nu}, \mathcal{Q} \mathcal{D}(u) \cdot \nu \big)_{|\Gamma_s}.
\end{align*}
\begin{theorem} \label{theo:stokes:2}
Let $\Omega$ be a bounded domain in $\R^n$, $n \ge 1$, with $C^{3}$-boundary
$\Gamma$ decomposing disjointly as $\Gamma = \Gamma_d \cup \Gamma_s$ with dist$\,(\Gamma_d,\Gamma_s)>0$.
Let $J=[0,T]$ and assume that $\rho_0$, $\eta_0$, $\lambda_0 \in
\cont{2}{}{\ol{\Omega}} \cap \sob{3}{2}{\Omega}$ as well as $\rho_0(x)$, $\eta_0(x)$,
$2 \eta_0(x) + \lambda_0(x) > 0$ for all $x \in \ol{\Omega}$.
Then the subproblem for $u$ in \eqref{eq:linear} possesses a unique solution in
\begin{align*}
Z_{1,\mathcal{B}}(J):= \big \{
v \in Z_1(J): &\;\; \mathcal{B}_d v \in
\sobb{\frac{5}{4}}{2}{J;\lp{2}{\Gamma_d;\R^n}}, \\
&\;\;\;\mathcal{B}_s v \in \sobb{\frac{5}{4}}{2}{J;\lp{2}{\Gamma_s}} \times
\sobb{ \frac{3}{4} }{2}{J;\lp{2}{\Gamma_s;\R^n}} \big \},\\
Z_1(J) = \sob{5/4}{2}{J;\lp{2}{\Omega;\R^n}} & \cap \sob{1}{2}{J;\sob{1}{2}{\Omega;\R^n}} \cap \lp{2}{J;\sob{3}{2}{\Omega;\R^n}},
\end{align*}
if and only if the data  $\tilde{f} := f - \mathcal{B}(D)c - \mathcal{C} \mu$,
$h_d$, $h_s=(h_{s1},h_{s2})$, $u_0$ satisfy the following conditions
\begin{compactenum}
\item $\tilde{f} \in \sob{1/4}{2}{J;\lp{2}{\Omega;\R^n}} \cap \lp{2}{J;\sob{1}{2}{\Omega;\R^n}}$;
\item $(h_d,h_{s1},h_{s2}) \in \mathbb{Y}_{0,d}(J;\R^n) \times
\mathbb{Y}_{0,s}(J) \times \mathbb{Y}_{1,s}(J;\R^n)$ with
$\mathbb{Y}_{k,i}(J;E):=$ \\
$\sobb{\frac{1}{2} ( 3 - k - \frac{1}{2})}{2}{J;\lp{2}{\Gamma_i;E}}
\cap \lp{2}{J;\sobb{3-k-\frac{1}{2}}{2}{\Gamma_i;E}}$, $k=0,1$, $i=d,s$, $E \in \{\R^n,\R\}$;
\item $u_0 \in \sobb{2}{2}{\Omega;\R^n}$;
\item $u_{0|\Gamma_d} = h_{d|t=0}$ in $\sobb{3/2}{2}{\Gamma_d;\R^n}$;
\item $\ska{ u_0 }{ \nu }{|\Gamma_s} = h_{s1|t=0}$ in $\sobb{3/2}{2}{\Gamma_s}$,
$\mathcal{Q} \mathcal{D}(u)_{|t=0} \cdot \nu_{|\Gamma_s} = h_{s2|t=0}$ in $\sobb{1/2}{2}{\Gamma_s;\R^n}$.
\end{compactenum}
\end{theorem}
\begin{theorem} \label{theo:stokes:3}
Let $\Omega$ be a bounded domain in $\R^n$, $n \ge 1$, with $C^{4}$-boundary
$\Gamma$ decomposing disjointly as $\Gamma = \Gamma_d \cup \Gamma_s$ with dist$\,(\Gamma_d,\Gamma_s)>0$.
Let $J=[0,T]$ and $p > \max \{ \frac{4}{3}, \frac{n}{3} \}$ with
$p \not = \frac{3}{2}$, $3$.
Further, assume that $\rho_0$, $\eta_0$, $\lambda_0 \in
\cont{2}{}{\ol{\Omega}} \cap \sob{3}{p}{\Omega}$ as well as $\rho_0(x)$, $\eta_0(x)$,
$2 \eta_0(x) + \lambda_0(x) > 0$ for all $x \in \ol{\Omega}$.
Then the subproblem for $u$ in \eqref{eq:linear} possesses a unique solution in
\begin{align*}
\mathcal{Z}_{1,\mathcal{B}}(J):= \big \{
v \in \mathcal{Z}_1(J): &\;\; \mathcal{B}_d v \in
\sobb{2-\frac{1}{2p}}{p}{J;\lp{p}{\Gamma_d;\R^n}}, \\
&\;\;\;\mathcal{B}_s v \in \sobb{2-\frac{1}{2p}}{p}{J;\lp{p}{\Gamma_s}} \times
\sobb{ \frac{3}{2}-\frac{1}{2p} }{p}{J;\lp{p}{\Gamma_s;\R^n}} \big \},
\end{align*}
if and only if the data  $\tilde{f} := f - \mathcal{B}(D)c - \mathcal{C} \mu$,
$h_d$, $h_s=(h_{s1},h_{s2})$, $u_0$ satisfy the following conditions
\begin{compactenum}
\item $\tilde{f} \in \mathcal{X}^n_{1,\Gamma}(J) := \{ \varphi \in \mathcal{X}_1^n(J): \,
\varphi_{|t=0,\Gamma_d} \in \sobb{2-3/p}{p}{\Gamma_d;\R^n}, \;
\ska{ \varphi_{|t=0} }{ \nu }{ |\Gamma_s } \in \sobb{2-3/p}{p}{\Gamma_s},
\\
\mathcal{Q} \mathcal{D}(\varphi)_{|t=0,\Gamma_s} \in \sobb{1-3/p}{p}{\Gamma_s;\R^n} \}$;
\item $(h_d,h_{s1},h_{s2}) \in \mathcal{Y}_{0,d}(J;\R^n) \times
\mathcal{Y}_{0,s}(J) \times \mathcal{Y}_{1,s}(J;\R^n)$ with
$\mathcal{Y}_{k,i}(J;E):=$ \\
$\sobb{2 - \frac{k}{2} - \frac{1}{2p}}{p}{J;\lp{p}{\Gamma_i;E}}
\cap \lp{p}{J;\sobb{4-k-\frac{1}{p}}{p}{\Gamma_i;E}}$, $k=0,1$, $i=d,s$, $E \in \{\R^n,\R\}$;
\item $u_0 \in \sobb{4-\frac{2}{p}}{p}{\Omega;\R^n}$;
\item $u_{0|\Gamma_d} = h_{d|t=0}$ in $\sobb{4-\frac{3}{p}}{p}{\Gamma_d;\R^n}$;
\item $\ska{  u_{0} }{ \nu }{ |\Gamma_s} = h_{s1|t=0}$ in $\sobb{4-\frac{3}{p}}{p}{\Gamma_s}$,
$\mathcal{Q} \mathcal{D}(u)_{|t=0} \cdot \nu_{|\Gamma_s} = h_{s2|t=0}$ in
$\sobb{3-\frac{3}{p}}{p}{\Gamma_s;\R^n}$;
\item $\rho_{0|\Gamma_d} \del{}{t} h_{d|t=0}
+ \mathcal{A}(D)u_{|t=0,\Gamma_d} = \tilde{f}_{|t=0,\Gamma_d}$
in $\sobb{2-\frac{3}{p}}{p}{\Gamma_d;\R^n}$ and
\\
$\rho_{0|\Gamma_s} \del{}{t} h_{s1|t=0} + \ska{ \mathcal{A}(D)u_{|t=0} }{\nu }{|\Gamma_s} =
\ska{ \tilde{f}_{|t=0} }{ \nu }{|\Gamma_s}$ in $\sobb{2 - \frac{3}{p}}{p}{\Gamma_s}$ if $p > \tfrac{3}{2}$;
\item $\del{}{t} h_{s2|t=0} + \mathcal{Q} \mathcal{D} ( \rho_0^{-1} \mathcal{A}(D)u )_{|t=0,\Gamma_s}
\cdot \nu_{|\Gamma_s} = \mathcal{Q} \mathcal{D} ( \rho_0^{-1} \tilde{ f} )_{|t=0,\Gamma_s}
\cdot \nu_{|\Gamma_s}$ in $\sobb{1 - \frac{3}{p}}{p}{\Gamma_s;\R^n}$ if $p > 3$.
\end{compactenum}
\end{theorem}
We will only give a proof for the last theorem. Theorem \ref{theo:stokes:2} can be proved by the same methods, 
the proof being even much simpler than that for the last theorem, due to the lower degree of regularity.

{\it Proof of Theorem \ref{theo:stokes:3}}. The necessity part is a consequence of trace theory and Sobolev embeddings.
The compatibility conditions follow by taking all possible temporal and spatial traces in
the differential equation as well as the boundary and initial conditions. Note that the higher
time regularity of $h_d,h_{s1},h_{s2}$ results from the additional time regularity of $\mathcal{B}_d v$
and $\mathcal{B}_s v$.

Concerning the sufficiency part, we know already from Theorem \ref{theo:stokes:1} that the subproblem
for $u$ in \eqref{eq:linear} possesses a unique solution $u$ in the regularity class $Z_1(J)$ so that it only remains to show that
$u$ enjoys the higher regularity properties stated in Theorem \ref{theo:stokes:3}. To achieve this,
one can use the well-known localization technique, flattening of the boundary, and perturbation arguments 
(see e.g. \cite{dhp-2}) to reduce the problem to associated full and half space problems with constant coefficients.
For these types of problems explicit solution formulas are available, which allow us to prove the desired extra regularity of $u$.
For the sake of brevity, we will only deal with the model problem in the half space
that comes from the boundary condition on $\Gamma_s$. The other half space case and the full space case are even simpler
and can be treated in the same way.

The problem in the half space $\R^n_+ := \R^{n-1} \times \R_+$ we have to study in order to be able to
treat the local problems involving the boundary segment $\Gamma_s$ reads as follows
\begin{equation} \label{eq:stokes}
\begin{aligned}
\del{}{t} u+u - \eta_0 \Delta u - (\lambda_0 + \eta_0)
\nabla \divv u & = f(t,x,y),
& \quad & t\in J, \; x \in \R^{n-1},\; y > 0,
\\
- \del{}{y} u^t = \theta(t,x), \quad  u^n & = \vartheta(t,x),
& \quad &  t\in J, \; x \in \R^{n-1},\; y = 0,
\\
u & = u_0(x,y),
& \quad & t = 0, \; x \in \R^{n-1},\; y>0,
\end{aligned}
\end{equation}
where we have set $u=(u^t,u^n) \in \R^{n-1} \times \R$. The outer unit normal in this case is given by $\nu = (0,\ldots,0,-1)^T$.
For the sake of convenience we added the zeroth order term $u$ on the left-hand side of the PDE, which is always possible, since
we have a compact time interval.
Transferring the regularity assumptions from the original problem to this
half space problem we obtain
\begin{align*}
& f \in \sob{1/2}{p}{J;\lp{p}{\R^n_+}} \cap \lp{p}{J;\sob{2}{p}{\R^n_+}},
f^n_{|t=0,y=0} \in \sobb{2-3/p}{p}{\R^n},
\\
& \del{}{y} f^t_{|t=0,y=0} \in \sobb{1-3/p}{p}{\R^n;\R^{n-1}},
\quad
u_0 \in \sobb{4-2/p}{p}{\R^n_+;\R^n},\\
& \vartheta \in \sobb{ (4-1/p) \frac{1}{2} }{p}{J;\lp{p}{\R^{n-1}}}
\cap \lp{p}{J;\sobb{4-1/p}{p}{\R^{n-1}}},
\quad
\del{\beta}{x} \vartheta \in  Y_0(J),
\\
& \theta \in \sobb{ (3-1/p) \frac{1}{2} }{p}{J;\lp{p}{\R^{n-1};\R^{n-1}}}
\cap \lp{p}{J;\sobb{3-1/p}{p}{\R^{n-1};\R^{n-1}}},
\quad
\del{\beta}{x} \theta \in Y_1(J;\R^{n-1}),
\\
& Y_i(J;E) := \sobb{ (2-i-1/p) \frac{1}{2} }{p}{J;\lp{p}{\R^{n-1};E}} \cap
\lp{p}{J;\sobb{2-i-1/p}{p}{\R^{n-1};E}},
\quad i=0,1,
\end{align*}
where $\del{\beta}{x}$ with $\beta \in \N^{n-1}$, $|\beta| \le 2$,
denotes tangential derivatives up to order $2$.
Moreover, the corresponding compatibility conditions take the form
\begin{align}
& u^n_{0|y=0} = \vartheta_{|t=0} \in \sobb{4-3/p}{p}{\R^n},
\quad
- [\del{}{y} u_0^t]_{|y=0} = \theta_{|t=0} \in \sobb{3-3/p}{p}{\R^n;\R^{n-1}},
\notag
\\
& [\del{}{t} \vartheta+\vartheta]_{|t=0} - \eta_0 [\Delta u^n_{0}]_{|y=0} =
[(\lambda_0 + \eta_0) \del{}{y} \divv u_0 + f^n_{|t=0}]_{|y=0}
\in \sobb{2-3/p}{p}{\R^n}\; \mbox{if}\;p>\frac{3}{2},
\label{cc:hs}
\\
& [\del{}{t} \theta+\theta]_{|t=0} + \eta_0 [ \del{}{y} \Delta u^t_0]_{|y=0} =
- [(\eta_0 + \lambda_0) \del{}{y} \nabla_x \divv u_{0}
+ \del{}{y} f^t_{|t=0}]_{|y=0} \in \sobb{1-3/p}{p}{\R^n;\R^{n-1}}
\notag
\end{align}
if $p>3$.

By maximal regularity with $\lp{p}{J;\lp{p}{\R^n_+;\R^n}}$ as base space for the PDE we know already that $u \in
\tilde{Z}(J) := \sob{ 1 }{p}{J;\lp{p}{\R^n_+;\R^n}} \cap \lp{p}{J;\sob{2}{p}{\R^n_+;\R^n}}$. Differentiating all equations of
\eqref{eq:stokes} with respect to the tangential variables, we see that $\del{\beta}{x}u \in \tilde{Z}(J)$ for all $|\beta| \le 2$,
by maximal $\lp{p}{J;\lp{p}{\R^n_+;\R^n}}$-regularity. Therefore
\begin{align*}
u \in \sob{1}{p}{J;\sob{2}{p}{\R^{n-1};\lp{p}{\R_+}}}
\cap
\lp{p}{J;\sob{2}{p}{\R^{n-1};\sob{2}{p}{\R_+}}}.
\end{align*}
Hence it remains to show that the normal derivatives $\del{j}{y} u$, $j \in\{1,2\}$,
lie in $\tilde{Z}(J)$ and that $\del{}{t} u \in \sob{1/2}{p}{J;\lp{p}{\R^n_+;\R^n}} \cap
\lp{p}{J;\sob{2}{p}{\R^n_+;\R^n}}$.
To establish this regularity, we will derive a solution formula for \eqref{eq:stokes}
from which the regularity can be read off.

To begin with, it is useful to consider
$v := \divv u$, which solves the problem
\begin{equation} \label{eq:heat}
\begin{aligned}
\del{}{t} v+v -  (2\eta_0 + \lambda_0) \Delta v &
= \nabla_x \cdot f^t + \del{}{y} f^n,
& \quad & t\in J, \; x \in \R^{n-1},\; y>0,
\\
- \del{}{y}v & = \psi,
& \quad & t\in J, \; x \in \R^{n-1},\; y=0,
\\
v & = \divv u_0 =: v_0,
& \quad & t = 0, \; x \in \R^{n-1},\; y>0,
\end{aligned}
\end{equation}
where
\begin{align*}
\psi & \,= - \nabla_x \cdot \del{}{y} u^t_{|y=0} - \del{2}{y} u^n_{|y=0}\\
 & \,= \nabla_x \cdot \theta + (2 \eta_0 + \lambda_0)^{-1}
[ f^n_{|y=0} - \del{}{t} \vartheta-\vartheta  + \eta_0 \Delta_x \vartheta -
( \eta_0 + \lambda_0 ) \nabla_x \cdot \theta ],
\end{align*}
in view of the identity
\begin{align*} - (2 \eta_0 + \lambda_0) \del{2}{y} u^n =  \eta_0\Delta_x u^n
+ (\eta_0 + \lambda_0) \nabla_x \cdot \del{}{y} u^t
- \del{}{t} u^n-u^n + f^n.
\end{align*}
Observe that $\psi$ belongs to $\sobb{1/2-1/4p}{p}{J;\lp{p}{\R^{n-1}}} \cap
\lp{p}{J;\sobb{2-1/p}{p}{\R^{n-1}}}$, which comes from $f^n_{|y=0}$ having the
least regularity. Further, the compatibility condition $-[\del{}{y} v_0]_{|y=0}
= \psi_{|t=0}$ is satisfied if $p>3/2$, by the first three conditions in \eqref{cc:hs}.
A solution formula for \eqref{eq:heat} is well-known, cf. \cite{pruess-1},
however we need a representation which allows to verify the desired
higher regularity and thus takes into account the corresponding
compatibility conditions. To be precise, the aim is to show the regularity
\begin{equation} \label{vreg}
\del{}{y} v,\; B^{1/2} v \in \sob{1/2}{p}{J;\lp{p}{\R^n_+}} \cap \lp{p}{J;\sob{2}{p}{\R^n_+}},
\end{equation}
where the operator $B$ is defined as $B=B_a = -\Delta_x +a^{-1}I $ with domain
$D(B)=\sob{2}{p}{\R^{n-1}}$ and $a=2 \eta_0 + \lambda_0$;
actually here we have to take the natural extension of
 $B$ to the corresponding space of functions that also depend on $t$ and $y$.
In what follows we will use the same symbol for both objects, whenever there is no danger of confusion.
Having established (\ref{vreg}) we will be able to infer the desired regularity of $u$
from the model problem for $u$ with right-hand side depending on $v$ which we consider below.

We introduce the operator $A_{\scriptscriptstyle 1/2}=\frac{1}{2}B-\partial_y^2$ with domain
$D(A_{\scriptscriptstyle 1/2} ) = \{ \varphi \in \sob{2}{p}{\R^n_+}:  \varphi_{|y=0} = 0\}$. 
Let $\phi$ denote the unique solution of
\begin{equation} \label{eq:phi}
\begin{split}
-\del{2}{y} \phi + \frac{1}{2} B \phi & =  e^{- (B/2)^{1/2}  y} g, \quad y > 0,
\quad
g := g(v_{0}) := [ \tfrac{1}{2} B v_0 - \del{2}{y} v_0]_{|y=0},
\\
\phi(0) & = v_{0|y=0},
\end{split}
\end{equation}
which is given by
\begin{align} \label{formula:phi}
\phi = \Phi(y)  v_{0|y=0} + \tfrac{y}{2} (\tfrac{1}{2} B)^{-1/2} \Phi(y) g
=
\Phi(y) v_{0|y=0} + A_{\scriptscriptstyle 1/2}^{-1} \Phi(y) g,
\end{align}
where $\Phi(y)$ denotes the analytic semigroup $e^{- (B/2)^{1/2} y}$. Then $\phi$
belongs to $\sobb{3-2/p}{p}{\R^n_+}$, since
$v_{0|y=0} \in \sobb{3-3/p}{p}{\R^{n-1}}$,
$g \in \sobb{1-3/p}{p}{\R^{n-1}}$, and by the mapping properties
of $\Phi$, see Proposition \ref{prop:1}. It follows from the construction of $\phi$ that
$v_0 - \phi \in A_{\scriptscriptstyle 1/2} ^{-1/2} D_{A_{\scriptscriptstyle 1/2} }(1-1/p,p)
= \{ \varphi \in \sobb{3-2/p}{p}{\R^n_+}:
\varphi_{|y=0} = A_{\scriptscriptstyle 1/2}  \varphi_{|y=0} = 0\}$, whenever traces make sense.

We next define for $\alpha>0$ the operators $S_{\alpha}(t) :=  e^{- \alpha(B/2) t}$,
and $T_\alpha(t) := e^{-\alpha A t}$, where we have set $A := A_1 = B -\del{2}{y}$
with domain $D(A) = D(A_{\scriptscriptstyle 1/2} )$. Let further $G = \del{}{t}$ with domain
$D(G)= \zero{ \sob{1}{p}{J;X} } := \{ \varphi \in \sob{1}{p}{J:X}: \varphi_{|t=0} = 0 \}$
and $F_\alpha := ( \alpha^{-1} G + B)^{1/2}$, $\alpha > 0$, with domain
$D(F_\alpha) = D(G^{1/2}) \cap D(B^{1/2}) =
\zero{ \sob{1/2}{p}{J;\lp{p}{\R^{n-1}}} } \cap \lp{p}{J;\sob{1}{p}{\R^{n-1}}}$.
These operators are sectorial, invertible and belong to $\bip{\lp{p}{J;\lp{p}{\R^{n-1}}}}$
with power angles $\theta_G \le \pi/2$ and $\theta_{F_\alpha} \le \pi/4$,
respectively; here we apply the Dore-Venni Theorem similarly as in \cite{pruess-1}
to show that $F_\alpha$ enjoys the claimed properties.
By means of these operators $v$ can be written as
\begin{align*}
v = \, &\, v_1 + e^{-F_a y} F_a^{-1} [ \psi - \del{}{y} v_{1|y=0} ],\\
 v_1  :=\, &\,
T_{a}(t) [ v_0 - \phi ] + S_{a}(t) \phi
\\
& +T_{a} * \left \{ \divv f  - e^{-F_{a} y} ( \divv f_{|y=0} -
S_{a}(\cdot) \divv f_{|y=0,t=0}) - S_{a}(\cdot) \Phi(y) \divv f_{|y=0,t=0} \right \}(t)
\\
& +t S_{a}(t) \Phi(y) \big[ \divv f_{|y=0,t=0} - \tfrac{1}{2} B  v_{0|y=0}
+ \del{2}{y} v_{0|y=0} \big] \\
& +\tfrac{y}{2} e^{-F_{a} y} F_{a}^{-1}
\big[ \divv f_{|y=0} - S_{a}(t) \divv f_{|y=0,t=0} \big].
\end{align*}

To see that $v$ possesses the regularity as mentioned, we remark that $\nabla \cdot f$
belongs to $\sob{1/4}{p}{J;\lp{p}{\R^n_+}} \cap \lp{p}{J;\sob{1}{p}{\R^n_+}}$
and thus, by trace theory, we obtain  $\nabla \cdot f_{|y=0} \in
\sobb{1/4-1/4p}{p}{J;\lp{p}{\R^{n-1}}} \cap \lp{p}{J;\sobb{1-1/p}{p}{\R^{n-1}}}$,
$\nabla \cdot f_{|t=0,y=0} \in \sobb{1-5/p}{p}{\R^{n-1}}$. The verification of the desired
regularity for $v$ is quite similar to the argument given below for $u^t$, which is why we omit the details here.

Relying on the regularity of $v$ just shown, it is now not so difficult to infer the claimed regularity properties of $u$. In fact,
$u$ can be considered as the solution of
\begin{alignat*}{3}
\del{}{t} u^t+u^t - \eta_0 \Delta u^t  &
= (\lambda_0 + \eta_0)  \nabla_x v + f^t(t,x,y) =: h^t(t,x,y),
& \quad & t > 0,\; x \in \R^{n-1},\;y>0,
\\
\del{}{t} u^n+u^n - \eta_0 \Delta u^n  &
= (\lambda_0 + \eta_0)  \del{}{y} v + f^n(t,x,y) =: h^n(t,x,y),
& \quad & t > 0,\;x \in \R^{n-1},\;y>0,
\\
- \del{}{y} u^t & = \theta(t,x), \quad  u^n = \vartheta(t,x),
& \quad & t > 0,\;x \in \R^{n-1},\;y=0,
\\
u^t & = u^t_0(y,x), \quad u^n = u^n_0(y,x)
& \quad & t = 0,\; x \in \R^{n-1},\;y>0.
\end{alignat*}
where we split again $u=(u^t,u^n)$ and $f=(f^t,f^n)$, and view
$\nabla v$, the regularity of which is known by means of the results above, as an inhomogeneity.
Therefore, the problem for $u^t$ and $u^n$ decouples and we have the representations
\begin{align*}
u^t = \,&\, u^t_1 + e^{-F_{\eta_0} y} F_{\eta_0}^{-1} \big[ \theta - \del{}{y} u^t_{1|y=0}\big ],\\
u^t_1 :=\,´&\, T_{\eta_0}(t) [ u^t_0 - \phi^t ]+ S_{\eta_0}(t) \phi^t
\\
&+T_{\eta_0} * \left \{
h^t - e^{-F_{\eta_0} y} (h^t_{|y=0} - S_{\eta_0}(\cdot) h^t_{|y=0,t=0})
- S_{\eta_0}(\cdot) \Phi(y) h^t_{|y=0,t=0}
\right \}(t)
\\
&+ t S_{\eta_0}(t) \Phi(y) \big[ h^t_{|y=0,t=0} - \tfrac{1}{2} B  u^t_{0|y=0}
+ \del{2}{y} u^t_{0|y=0} \big]\\
& + \tfrac{y}{2} e^{-F_{\eta_0} y} F_{\eta_0}^{-1}
\big[ h^t_{|y=0} - S_{\eta_0}(t) h^t_{|y=0,t=0} \big],
\end{align*}
and
\begin{align*}
u^n =\, &\, u^n_1 + T_{\eta_0} * \left \{
h^n - e^{-F_{\eta_0} y} \big(h^n_{|y=0} - S_{\eta_0}(\cdot) h^n_{|y=0,t=0}\big) - S_{\eta_0}(\cdot)
\Phi(y) h^n_{|y=0,t=0}
\right \} (t)
\\
&+ T_{\eta_0}(t) [ u^n_0 - \phi^n ]
+ \tfrac{y}{2} e^{-F_{\eta_0} y} F_{\eta_0}^{-1}
\big[ h^n_{|y=0} - S_{\eta_0}(t) h^n_{|y=0,t=0}\big]
+ e^{-F_{\eta_0} y} [ \vartheta - u^n_{1|y=0} ] ,
\\
u^n_1 :=\,&\, S_{\eta_0}(t) \phi^n +
t S_{\eta_0}(t) \Phi(y) \big[ h^n_{|y=0,t=0} -\frac{1}{2} B u^n_{0|y=0} + \del{2}{y} u^n_{0|y=0} \big].
\end{align*}
Here, $B=B_{\eta_0}$ and $\phi = (\phi^t,\phi^n) \in \sobb{4-2/p}{p}{\R^n_+;\R^{n-1}} \times \sobb{4-2/p}{p}{\R^n_+}$
denotes the solution of \eqref{eq:phi} with the data given there replaced by $(\Phi(y) g(u_0), u_{0|y=0})$, which
implies $u_0 - \phi \in A^{-1} D_{A}(1-1/p,p) =
\{ \varphi \in \sobb{4-2/p}{p}{\R^n_+;\R^n}: \varphi_{|y=0} = A \varphi_{|y=0} = 0 \}$.
To understand where this regularity comes from,
one first verifies that $u_{0|y=0} \in \sobb{4-3/p}{p}{\R^{n-1};\R^n}$,
$g(u_0) \in \sobb{2-3/p}{p}{\R^{n-1};\R^n}$, and $\Phi(\cdot) g(u_0) \in
\sobb{2-3/p}{p}{\R^n_+;\R^n}$, by applying Proposition \ref{prop:1}. Further, due to the second
representation of $\phi$, see \eqref{formula:phi}, taking the derivative $\del{}{y}$ corresponds to applying
$B^{1/2}$ in terms of regularity, and therefore $B \phi$, $\del{2}{y} \phi \sim
\Phi(y) B u_{0|y=0} + B A_{1/2}^{-1} \Phi(y) g$,
where both terms lie in $\sobb{2-2/p}{p}{\R^n_+;\R^n}$.

Observe that the compatibility conditions were incorporated in the above solution formulas to the result that
$\theta_{|t=0} = \del{}{y} u_{1|y=0,t=0}^t$, $\del{}{t} \theta_{|t=0} = \del{}{t}
\del{}{y} u_{1|y=0,t=0}^t$, $\vartheta_{|t=0} = u^n_{1|y=0,t=0}$, and
$\del{}{t} \vartheta_{|t=0} = \del{}{t} u^n_{1|y=0,t=0}$.

The solution formulas allow us to verify the desired additional regularity of $u$, that is
$\del{k}{x_i} u$,
$\del{k}{y} u \in \tilde{Z}(J)$, $k=1,2$, and $\del{}{t} u \in \sob{1/2}{p}{J;\lp{p}{\R^n_+;\R^n}}$.
We will give the argument for the component $u^t$.
The component $u^n$ can be discussed similarly. In what follows we will also
use the symbol $\tilde{Z}(J)$ for functions taking values in $\R^{n-1}$.

We first study the term
$w_1 := T_{\eta_0}(t)[u^t_0 - \phi^t]$. Here taking the derivatives $\del{}{t}$, $\del{2}{x_i}$,
or $\del{2}{y}$ corresponds to applying $A$. Therefore, from $A [u^t_0 - \phi^t] \in D_A(1-1/p,p)$
we conclude that $w_1 \in \sob{2}{p}{J;\lp{p}{\R^n_+;\R^{n-1}}} \cap
\lp{p}{J;\sob{4}{p}{\R^n_+;\R^{n-1}}}$.

The third term in the representation of $u^t_1$, namely
$w_3:=T_{\eta_0} * \{ \ldots \}$,
is more involved. To begin with, observe that $\{ \ldots \}_{|y=0} = 0$, which along with
$\{ \ldots \} \in \lp{p}{J;\sob{2}{p}{\R^n_+;\R^{n-1}}}$ leads to
$A w_3 \in \tilde{Z}(J)$. Moreover, it is easy to see that
$S_{\eta_0}(t) \Phi(y) h^t_{|y=0,t=0} \in \tilde{Z}(J)$,
$e^{-F_{\eta_0} y} (h^t_{|y=0} - S_{\eta_0}(t) h^t_{|y=0,t=0}) \in
\sobb{1/2+1/2p}{p}{J;\lp{p}{\R^n_+;\R^{n-1}}} \cap \lp{p}{J;\sob{2}{p}{\R^n_+;\R^{n-1}}}$,
and since $h^t$ belongs to $\sob{1/2}{p}{J;\lp{p}{\R^n_+;\R^{n-1}}} \cap
\lp{p}{J;\sob{2}{p}{\R^n_+;\R^{n-1}}}$, we may conclude from the identity
$\del{}{t} w_2 = T * A \{ \ldots \} + \{ \ldots \}$ the desired regularity.

The next term we study is $w_2 :=S_{\eta_0}(t) \phi^t$, it even lies in
$\tilde{Z}^2(J) := \sob{2}{p}{J;\lp{p}{\R^n_+;\R^{n-1}}} \cap \lp{p}{J;\sob{4}{p}{\R^n_+;\R^{n-1}}}$.
To see this, observe first that we have $\phi^t \in \sobb{4-2/p}{p}{\R^n_+;\R^{n-1}}$ and
$S_{\eta_0}(t) [B \phi^t]_{|y=0} \in \lp{p}{J;\sobb{2-1/p}{p}{\R^n;\R^{n-1}}} \cap
\sobb{1-1/2p}{p}{J;\lp{p}{\R^n;\R^{n-1}}}$, by Proposition \ref{prop:1}.
Further, as we have the relations $\del{}{t} w_2 \sim B w_2$,
$\del{2}{x_i} w_2 \sim B w_2$, and $\del{2}{y} w_2 \sim B w_2$, it is sufficient
to consider $\hat{v} := B w_2$, which solves
\begin{align*}
\del{}{t}\hat{ v}+\hat{v} - \eta_0 \Delta \hat{v}  & = 0,
\quad t > 0,\; x \in \R^{n-1},\;y>0,
\\
\hat{v}_{|y=0} & \in \sobb{1-1/2p}{p}{J;\lp{p}{\R^n;\R^{n-1}}} \cap
\lp{p}{J;\sobb{2-1/p}{p}{\R^n;\R^{n-1}}},
\\
\hat{v}_{|t=0} &  \in \sobb{2-2/p}{p}{\R^n_+;\R^{n-1}}.
\end{align*}
By maximal $L_p$-regularity it follows that $\hat{v} \in \tilde{Z}(J)$, and hence $w_2 \in \tilde{Z}^2(J)$.

We next investigate the term
$w_4  := t S_{\eta_0}(t) \Phi(y) [ h^t_{|y=0,t=0} - \tfrac{1}{2} B  u^t_{0|y=0}
+ \del{2}{y} u^t_{0|y=0} ] =: t S_{\eta_0}(t) \Phi(y) \tilde{w}_{00} =: t \tilde{w}$. It belongs
to $\tilde{Z}^2(J)$ as well. In fact, observe that
$\tilde{w}_{00}$ lies in $\sobb{2-3/p}{p}{\R^n;\R^{n-1}}$ and $\tilde{w} \in Z(J)$, as $\tilde{w}$
solves the problem above with right side $0$, boundary data $S(t) \tilde{w}_{00}$, and
initial data $\Phi(y) \tilde{w}_{00}$, where the data possess the regularity stated above.
Moreover,  $w_4$ solves the problem above with right side $\tilde{w}$,  boundary data
$t S(t) \tilde{w}_{00}$, and initial data $0$. Because of this observation we are able to
rewrite $w_4$ as follows.
\begin{align*}
w_4(t,y) & = e^{-F_{\eta_0} y} t S_{\eta_0}(t) \tilde{w}_{00}
+ \tfrac{1}{2} F^{-1/2}_{\eta_0} \int_0^\infty
[ e^{-F_{\eta_0} |y-s| }- e^{-F_{\eta_0} (y+s)} ] \tilde{w}(t,s) \, ds
\\
& \equiv e^{-F_{\eta_0} y} (G + \tfrac{1}{2} B)^{-1} S_{\eta_0}(t) \tilde{w}_{00} + \ldots,
\end{align*}
and this representation reveals the claimed regularity.

Finally, we have to look at
$w_5(t,y) \! := \tfrac{y}{2} F_{\eta_0}^{-1} e^{-F_{\eta_0} y} [h^t_{|y=0} -
S_{\eta_0}(t) h^t_{|t=0,y=0}]$. Having in mind that $[h^t_{|y=0} - S_{\eta_0}(t) h^t_{|t=0,y=0}]
\in \zero{\sobb{1/2-1/4p}{p}{J;\lp{p}{\R^n;\R^{n-1}}}} \cap \lp{p}{J;\sobb{2-1/p}{p}{\R^n;\R^{n-1}}}$
and thus
\begin{align*}
e^{-F_{\eta_0} y} [ \ldots ] \in \zero{\sobb{1/2+1/2p}{p}{J;\lp{p}{\R^n_+;\R^{n-1}}}}
\cap \lp{p}{J;\sob{2}{p}{\R^n_+;\R^{n-1}}}
\end{align*}
as well as the embedding
\begin{align*}
\zero{\sobb{1/2+1/2p}{p}{J;\lp{p}{\R^n_+;\R^{n-1}}}}
\hookrightarrow
\zero{\sob{1/2}{p}{J;\lp{p}{\R^n_+;\R^{n-1}}}} ,
\end{align*}
one easily verifies that $\del{}{t} w_5$, $B w_5$, $\del{2}{y} w_5 \in
\sob{1/2}{p}{J;\lp{p}{\R^n_+;\R^{n-1}}} \cap \lp{p}{J;\sob{2}{p}{\R^n_+;\R^{n-1}}}$
which shows
\begin{align*}
w \in \sob{3/2}{p}{J;\lp{p}{\R^n_+;\R^{n-1}}} \cap
\sob{1}{p}{J;\sob{2}{p}{\R^n_+;\R^{n-1}}} \cap \lp{p}{J;\sob{4}{p}{\R^n_+;\R^{n-1}}}.
\end{align*}
This completes the proof. \eproof

\section{The continuity equation} \label{sec:CE}

In this section the equation of conservation of mass is carefully studied concerning the regularity dependence on $u$, 
that is we are interested in the optimal regularity of the density $\rho$ we can expect, given a certain regularity of the 
velocity $u$. Since a third order term of $\rho$ appears in the Cahn-Hilliard
equation, which is supposed to be in $\mathcal{X}_2(J)=\lp{p}{J;\lp{p}{\Omega}}$, we need
$\del{}{x_i} \del{}{x_j} \del{}{x_k} \rho \in \lp{p}{J;\lp{p}{\Omega}}$, $i,j,k \le 3$, at least.
A similar situa\-tion occurs in the Navier-Stokes equation, since here first order terms
of $\rho$ have to be in $\lp{p}{J;\sob{2}{p}{\Omega}}$. We will see that this
regularity and even more can be gained from the assumption $u \in
\mathcal{Z}_1(J)$.

\begin{lemma} \label{lem:rho}
Let $\Omega$ be a bounded domain in $\R^n$, with $C^4$ boundary $\Gamma$, $J_0=[0,T_0]$ with $T_0>0$ being fixed or $J_0=\R_+$,
$J=[0,T]\subset J_0$ a compact time interval, and $\hat{p} < p < \infty$. Let $\ol{u}\in\mathcal{Z}_1(J_0)$ be a fixed function. 
Further, assume that $\rho_0 \in \sob{3}{p}{\Omega}$ with $\rho_0(x) > 0$
for all $x \in \ol{\Omega}$, and $u \in \mathcal{Z}_1(J)$ satisfies
$\norm{u - \ol{u}}{ \zero{\mathcal{Z}_1(J)} } \le 1$ as well as $\ska{u}{\nu}{} \ge 0$ on $\Gamma$.
Then problem \eqref{eq:me} supplemented with initial condition $\rho(0) = \rho_0$
possesses a unique non-negative solution $\rho \in \mathcal{Z}_3(J)$, which defines
a linear solution operator $L[u]$ according to $\rho(t) = L[u](t) \rho_0$, and there exists
a constant $\varsigma > 0$ independent of $T$ and $u$ such that
\begin{align} \label{est:rho}
\norm{\rho}{ \mathcal{Z}_3(J) } \le \varsigma.
\end{align}
\end{lemma}
Before we give the proof we make an important comment on how to get estimates that are independent of $T$. 
This will be crucial also for the following sections. Suppose that $J$ and $J_0$ are as in Lemma \ref{lem:rho}.
Let $\ol{w} = (\ol{u},\ol{c},\ol{\mu}) \in
\mathcal{Z}_1(J_0) \times \mathcal{Z}_2(J_0) \times \mathcal{Z}_\mu(J_0)$ be a triple of fixed functions and consider 
an arbitrary triple $w=(u,c,\mu)\in \mathcal{Z}_1(J) \times \mathcal{Z}_2(J) \times \mathcal{Z}_\mu(J)$ satisfying 
$\norm{w - \ol{w}}{ \zero{\mathcal{Z}_1}(J) \times \zero{\mathcal{Z}_2(J) \times {}_0\mathcal{Z}_\mu(J)} } \le 1$, 
that is we have in particular $(u,\partial_t u, c,\mu)(0)-(\ol{u},\partial_t \ol{u},\ol{c},\ol{\mu})(0)  = 0$.
Set $\norm{\ol{w}}{\mathcal{Z}_1(J_0) \times \mathcal{Z}_2(J_0) \times \mathcal{Z}_\mu(J_0)} =: \omega$. 
Let $Y(J)$ denote a function space with the property $\mathcal{Z}_1(J) \times \mathcal{Z}_2(J)
\times \mathcal{Z}_\mu(J)\hookrightarrow Y(J)$ and assume that $\norm{{v}}{Y(J)}\le \norm{{v}}{Y(J_0)}$ for all $v\in Y(J_0)$. Let us further assume that
there exists a bounded extension operator $E_+$ from $\zero{ Y(J) }$ to $\zero{ Y }(J_0)$ satisfying
$\|E_+\|_{ \bdd{ \zero{\mathcal{Z}_1}(J_0) \times \zero{\mathcal{Z}_2(J_0) \times {}_0\mathcal{Z}_\mu(J_0)} }{
\zero{\mathcal{Z}_1}(J) \times \zero{\mathcal{Z}_2(J) \times {}_0\mathcal{Z}_\mu(J)} } } =: m_+ < \infty$, 
where $m_+$ does not depend on $T$. These assumptions make it possible to estimate as follows
\begin{equation} \label{est:uniform}
\begin{split}
\norm{w}{Y(J)} & \le \norm{w - \ol{w}}{ \zero{Y(J)} } +
\norm{\ol{w}}{Y(J_0)} \le \norm{E_+(w-\ol{w})}{ \zero{Y(J_0)} } + C_E \omega
\\
& \le C_E \norm{E_+(w-\ol{w})}{\zero{\mathcal{Z}_1(J_0)} \times \zero{\mathcal{Z}_2(J_0) \times {}_0\mathcal{Z}_\mu(J_0)}}
+ C_E \omega
\\
& \le C_E m_+ \norm{w-\ol{w}}{ \zero{ \mathcal{Z}_1(J)} \times \zero{\mathcal{Z}_2(J) \times {}_0\mathcal{Z}_\mu(J)} }
+ C_E \omega
\le C_E m_+ + C_E \omega =: k < \infty.
\end{split}
\end{equation}
In particular this applies to functions $w \in \Sigma_T$, and it is also clear that in this way we get a uniform 
estimate for $u$ from Lemma \ref{lem:rho}. In fact, the function spaces under consideration satisfy the
above conditions; corresponding extension operators with uniform bound w.r.t. $T$ can be constructed for
these functions spaces (with vanishing traces at $t=0$) by means of standard reflection techniques.

${}$

{\it Proof of Lemma \ref{lem:rho}.}  Existence and uniqueness of $\rho \in \cont{1}{}{J;\lp{p}{\Omega}} \cap
\cont{}{}{J;\sob{3}{p}{\Omega}}$ as well as the independence  of the corresponding estimate of $\rho$
on $T$ and $u$ has been shown in \cite[Lemma 4.1]{k1}.
Hence we only have to verify the additional regularity $\rho \in
\sob{2+1/4}{p}{J;\lp{p}{\Omega}} \cap \cont{1}{}{J;\sob{2}{p}{\Omega}}$.

{\it Step I - $\rho \in \cont{1}{}{J;\sob{2}{p}{\Omega}}$}
This regularity just follows from $\rho \in \cont{}{}{J;\sob{3}{p}{\Omega}}$
and $u \in \mathcal{Z}_1(J)$ and the equation
\begin{align*}
\del{}{t} \del{}{x_i} \del{}{x_j} \rho
& =
- \del{}{x_i} \del{}{x_j} \del{}{x_k} \rho u_k - \del{}{x_j} \del{}{x_k} \rho \del{}{x_i} u_k
- \del{}{x_i} \del{}{x_k} \rho \del{}{x_j} u_k - \del{}{x_k} \rho   \del{}{x_i} \del{}{x_j} u_k
\\
& \quad - \del{}{x_i} \del{}{x_j} \rho \del{}{x_k} u_k - \del{}{x_j} \rho \del{}{x_i} \del{}{x_k} u_k
- \del{}{x_i} \rho \del{}{x_j} \del{}{x_k} u_k  - \rho \del{}{x_i} \del{}{x_j} \del{}{x_k} u_k,
\end{align*}
where we made use of Einstein's summation convention.

{\it Step II - $\rho \in \sob{2+\frac{1}{4}}{p}{J;\lp{p}{\Omega}}$.}
This regularity follows from $\nabla \cdot \del{}{t}u \in
\sob{1/4}{p}{J;\lp{p}{\Omega}}$, which is ensured by
the embedding $\mathcal{ Z}_1(J) \hookrightarrow
\sob{1+1/4}{p}{J;\sob{1}{p}{\Omega;\R^n}}$,
and the equation
\begin{align} \label{eq:rho-tt}
\del{2}{t} \rho = - \nabla \del{}{t} \rho \cdot u  - \rho \divv \del{}{t} u
- \del{}{t} \rho \divv u - \nabla \rho \cdot \del{}{t} u,
\end{align}
by a bootstrap argument.
In fact, the idea is to study the regularity of the right-hand side of (\ref{eq:rho-tt}) which is
at most the one mentioned above. Since we do not yet know $\nabla \del{}{t} \rho \in \sob{1/4}{p}{J;\lp{p}{\Omega}}$
appearing on the right side of \eqref{eq:rho-tt}, we first look at the equation
\begin{align*}
\del{}{x_i} \del{}{t}  \rho
& =
- \del{}{x_i} \del{}{x_k} \rho u_k - \del{}{x_k} \rho  \del{}{x_i} u_k
- \del{}{x_i} \rho  \del{}{x_k} u_k - \rho  \del{}{x_i} \del{}{x_k} u_k
\end{align*}
in $\sob{1}{p}{J;\lp{p}{\Omega}}$.
Since $\rho \in \cont{1}{}{J;\lp{p}{\Omega}}$ and  $u \in \mathcal{Z}_1(J) \subset \sob{1}{p}{J;\sob{2}{p}{\Omega;\R^n}}$,
one easily infers that $\rho \in \sob{2}{p}{J;\sob{1}{p}{\Omega}}$. To show the regularity $\rho \in \sob{2+1/4}{p}{J;\lp{p}{\Omega}}$, 
we argue again with equation \eqref{eq:rho-tt}, but now using the newly gained regularity for $\rho$, together with the regularity of $u$.
\eproof

${}$

The next lemma concerns the estimate of differences of solutions
of the equation of mass. The proof is analogous to \cite[Lemma 4.3]{k1};
recall that
\begin{align*}
Z_1(J) & = \sob{5/4}{2}{J;\lp{2}{\Omega;\R^n}} \cap \sob{1}{2}{J;\sob{1}{2}{\Omega;\R^n}}\cap \lp{p}{J;\sob{3}{2}{\Omega;\R^n}},
\\
Z_3(J) & = \sob{2}{2}{J;\lp{2}{\Omega}} \cap \cont{1}{}{J;\sob{1}{2}{\Omega}} \cap \cont{}{}{J;\sob{2}{2}{\Omega}}.
\end{align*}
\begin{lemma} \label{lem:diff-rho}
Let $\Omega$, $J$, $J_0$, $\rho_0$, and $\ol{u}$ be as in the previous lemma.
Suppose that $u_1$, $u_2 \in \mathcal{Z}_1(J)$
with $\norm{u_i - \ol{u}}{ \zero{\mathcal{Z}_1(J)} } \le 1$ and
$\ska{u_i}{\nu}{} \ge 0$ on $\Gamma$, $i=1,2$. Further let
$\rho_i = L[u_i]\rho_o\in  \mathcal{Z}_3(J)$, $i=1,2$.
Then there is a constant $\kappa_1(T) > 0$ with the property
$\kappa_1(T) \to 0$ as $T \to 0$, such that
\begin{equation}  \label{est:diff-rho}
\norm{\rho_1 - \rho_2}{ \zero{Z_3(J)} } \le
\kappa_1(T) \norm{u_1 - u_2}{ \zero{Z_1(J)} }.
\end{equation}
\end{lemma}
{\it Proof}. Suppose that $(u_i,\rho_i) \in \mathcal{Z}_1(J) \times \mathcal{Z}_3(J)$
solve the equation of conservation of mass. Letting
$\varrho:= \rho_1 - \rho_2$ and $v:=u_1 - u_2$,  $(\varrho,v)$ satisfies
\begin{equation} \label{eq:diff-rho}
\begin{split}
\del{}{t} \varrho + \divv (\varrho u_1) & = - \divv( \rho_2 v),
\quad (t,x) \in J \times \Omega,
\\
\varrho|_{t=0} & = 0, \qquad \qquad \:\;\, x \in \Omega.
\end{split}
\end{equation}
Observe that the right-hand side $\divv (\rho_2 v)$ belongs to
$\lp{2}{J;\sob{2}{2}{\Omega}}$, since the difference $v$ is considered in $\zero{Z_1(J)}$.
To establish the estimate \eqref{est:diff-rho}, one only has to
adopt the proof of \cite[Lemma 4.3]{k1} by using \cite[Lemma 4.2]{k1}. Time regularity, i.e.
$\del{}{t} \varrho \in \cont{}{}{J;\sob{1}{2}{\Omega}} \cap \sob{1}{2}{J;\lp{2}{\Omega}}$, follows
from the equation directly. \eproof
\section{An estimate in a weaker norm for Cahn-Hilliard} \label{sec:WCH}
In this section we will derive an estimate in $Z_2\times Z_\mu$ for differences of solutions $(c,\mu)$ to the
Cahn-Hilliard problem in \eqref{eq:nsch-2}. This estimate will be crucial for proving
the strict contraction property of the fixed point mapping $\mathcal{G}$ which was introduced in Section \ref{sec:3}.

Let $w_1,w_2\in \Sigma_T$, $\rho_i=L[u_i]\rho_0$, and $w_i'=\mathcal{G}(w_i)$,
$i=1,2$. The Cahn-Hilliard subproblem for $(c_i',\mu_i')$ reads as
\begin{alignat*}{2}
\tfrac{\varepsilon_0 \rho_0}{\gamma_0} \,\del{}{t} c_i'
-
\divv( \varepsilon_0 \nabla\mu_i' ) & = F_2(w_i,\rho_i),
 \quad &(t,x) \in J \times \Omega,
\\
-\mu_i'-\divv(\varepsilon_0 \nabla c_i') & = F_\mu(w_i,\rho_i), \quad &(t,x) \in J \times \Omega,\\
\partial_\nu c_i'=\partial_\nu \mu_i' & = 0,  \quad &(t,x)\in J\times \Gamma,\\
c_i'|_{t=0} & = c_0,  \quad & x \in \Omega.
\end{alignat*}
Next we set $c' := c_1'-c_2'$, $\mu' := \mu_1' - \mu_2'$, $u := u_1 - u_2$, $c:= c_1-c_2$, $\rho := \rho_1 - \rho_2$,
and $\mu :=\mu_1 - \mu_2$, and define the operators $M \varphi := \frac{\varepsilon_0 \rho_0}{\gamma_0} \varphi$ with domain
$D(M):= Y := \lp{2}{J;X}$ and $X:=\lp{2}{\Omega}$, $G \varphi:= \del{}{t} \varphi$ with natural domain
$D(G) := \zero{\sob{1}{2}{J;X}}$, and $\mathcal{A}_\alpha \varphi := - \divv \, ( \alpha \nabla \varphi)$,
$\alpha \in \{ \varepsilon_0, \gamma_0 - \gamma_1, \gamma_1 - \gamma_2, ...\}$
with domain $D(\mathcal{A}_\alpha) := \{ v \in \sob{2}{2}{\Omega}: \del{}{\nu} v = 0 \}$.
Further, let $A_\alpha$ denote the natural extension of  $\mathcal{A}_\alpha$ to $\lp{2}{J;X}$ with domain
$D(A_\alpha) = \lp{2}{J;D(\mathcal{A}_\alpha)}$. Then $M$, $G$, $A_{\varepsilon_0}+I$ are sectorial and invertible operators,
and $M$, $A_{\varepsilon_0}+I$ are self-adjoint. Moreover, these operators belong to the class
$\bip{Y}$ (which coincides with $\hinfty{ Y }$, since $Y$ is a Hilbert space)
with power angles $\theta_M = 0$, $\theta_G = \pi/2$, and $\theta_{A_{\varepsilon_0}+I} = 0$, respectively.

Inserting the second PDE into the first one and taking differences we obtain the following problem for $c'$.
\begin{equation} \label{eq:ch-weak-1}
\begin{split}
M G c' + A_{\varepsilon_0} [A_{\varepsilon_0} c'-\Phi_2 ]
& = \Phi_1
\end{split}
\end{equation}
where we set
\begin{align*}
\Phi_1 & := \tfrac{\varepsilon_0}{\gamma_0} \big ( G \phi_1 - \divv \phi_2 + A_{\gamma_0 - \gamma_1} \mu
- A_{\gamma_1 - \gamma_2}\mu_2 \big)
+ \tfrac{\varepsilon_0^2}{\gamma_0} \nabla ( \tfrac{\gamma_0}{\varepsilon_0} ) \cdot
\nabla \mu,
\\
\phi_1 & := (\rho_0 - \rho_1) c - \rho c_2, \quad \phi_2:= \rho_1 u_1 c + c_2 u_1 \rho + c_2 \rho_2 u,
\\
\Phi_2 & := F_\mu(w_1,\rho_1) - F_\mu(w_2,\rho_2) \equiv
- A_{\varepsilon_1 - \varepsilon_0} c - A_{\varepsilon_1 - \varepsilon_2} c_2
+
\rho^{-1}_1 \varepsilon_1 \nabla c_1 \cdot \nabla \rho
\\
& \phantom{:= \,}
+
[ \rho^{-1}_1 \varepsilon_1 \nabla c_1  - \rho^{-1}_2 \varepsilon_2 \nabla c_2] \cdot \nabla \rho_2
-
[ \del{}{c} \psi(\rho_1,c_1) -  \del{}{c} \psi(\rho_2,c_2) ].
\end{align*}
Once we have an estimate for the difference $c'$ in ${Z}_2(J)$, we obtain an estimate for the difference
$\mu'$ in $Z_\mu(J)$ by using
\begin{align} \label{eq:ch-weak-1a}
\mu' = A_{\varepsilon_0}c' - \Phi_2.
\end{align}
We have the following result.
\begin{lemma} \label{lem:diff-ch}
Let $\Omega \subset \R^n$, $n \ge 1$,  be a bounded domain with $C^4$
smooth boundary $\Gamma$, $J=[0,T]$ a compact time interval,
and $p \in (\hat{p},\infty)$. Assume that
\begin{itemize}
\item[(i)] $\rho_0 \in \sob{3}{p}{\Omega}$, $c_0 \in \sobb{4-4/p}{p}{\Omega}$;
\item[(ii)] $\varepsilon \in \cont{4}{}{\R^2}$, $\gamma \in \cont{2}{}{\R^2}$, $\ol{\psi}\in \cont{5}{}{\R^2}$;
\item[(iii)] $(u_i,c_i,\mu_i) \in \Sigma_T$, $\rho_i = L[u_i]\rho_0 \in \mathcal{Z}_3(J)$, $i=1,2$,
and the pair $(c',\mu') \in \zero{\mathcal{Z}_2(J)} \times \zero{\mathcal{Z}_\mu(J)}$ solves \eqref{eq:ch-weak-1}, \eqref{eq:ch-weak-1a}.
\end{itemize}
Then there exists a constant $\kappa_2(T) > 0$ with $\kappa_2(T) \to 0$ as $T
\to 0$, such that
\begin{equation} \label{est:diff-ch}
\norm{ (c',\mu') }{ \zero{Z_2}(J) \times Z_\mu(J) } \le \kappa_2(T)
\norm{(u,c,\mu)}{ \zero{ Z_1(J) } \times \zero{Z_2(J)} \times Z_\mu(J)} .
\end{equation}
\end{lemma}
{\it Proof.} {\bf a)} Let $G$ and $A_{\varepsilon_0}$ be defined as above.
Then $G^{3/4}$ is also a sectorial operator, it is invertible and belongs to $\mathcal{BIP}(Y)$ with
$\theta_{G^{3/4}} \le 3\pi/8$, and $D(G^{3/4})=\zero{ \sob{3/4}{2}{J;\lp{2}{\Omega}}}$. Moreover,
$( A_{\varepsilon_0} +I )^{3/2}$ is sectorial, belongs to $\mathcal{BIP}(Y)$ with $\theta_{ (A_{\varepsilon_0} + I )^{3/2} }=0$,
and $(A_{\varepsilon_0} + I)^{3/2} = (A_{\varepsilon_0} + I)^{1/2} (A_{\varepsilon_0} + I)$.
Since $G^{3/4}$ and $(A_{\varepsilon_0} + I)$ commute, the Dore-Venni Theorem (see
\cite{pruess-1}) implies that $B := G^{3/4} + (A_{\varepsilon_0} + I)^{3/2}$ with domain
\begin{align*}
D(B) = D(G^{3/4}) \cap D( (A_{\varepsilon_0} + I)^{3/2} )
& = \zero{ \sob{3/4}{2}{J;\lp{2}{\Omega}} } \cap \lp{2}{J;D( (\mathcal{A}_{\varepsilon_0}+I)^{3/2}) },
\\
D( (\mathcal{A}_{\varepsilon_0} +I)^{3/2} ) & = \{ \phi \in \sob{3}{2}{\Omega}: \del{}{\nu} \phi = 0 \}
\end{align*}
is invertible, sectorial, and belongs to $\mathcal{BIP}(Y)$ with angle
$\theta_B \le \max \{\theta_{G^{3/4}}, \theta_{ (A_{\varepsilon_0} + I)^{3/2}} \} = 3 \pi/8$. Let further
$F:=MG + A_{\varepsilon_0}^2$ with domain $D(F) = D(G) \cap D((A_{\varepsilon_0} + I)^2)$.
Maximal $L_p$-regularity of the Cahn-Hilliard problem then implies that $F \in \liss{ D(F) }{ Y }$.

A crucial point in the following argument will be that all subsequent constants $C$, $C_i$ etc.,
which may differ from line to line, are always independent of the unknowns and $T$.
This is possible due to working in function spaces with time trace zero at $t=0$.
(The point is that functions lying in such spaces can be extended to $\R_+$, where the extension
operator is bounded by a constant independent of $T$.)

{\bf b)} We next derive an estimate for $c'$ in $\zero{Z}_{2}(J)$ which is in terms of the operators just defined
(see \eqref{eq:norm-c-1} below), and thus more appropriate for inferring the desired estimate for $c'$
from problem \eqref{eq:ch-weak-1}. Note first that $\Phi_2 \in
\zero{\sob{1/2}{p}{J;\lp{p}{\Omega}} \cap \lp{p}{J;\sob{2}{p}{\Omega}}}$, $c' \in \zero{\mathcal{ Z}_2(J)}$, and
$\del{}{\nu} (A_{\varepsilon_0} c' - \Phi_2) = 0$ at $J\times \Gamma$, see also the estimates below.
From  \eqref{eq:ch-weak-1} it follows that $c' - (A_{\varepsilon_0} +I)^{1/2} B^{-1} \Phi_2 \in D(F)$; in fact we have
$(A_{\varepsilon_0} +I)^{1/2} B^{-1} \Phi_2 \in \zero{\mathcal{Z}_2(J)}$ as well as $\del{}{\nu} c'  = 0$,
$\del{}{\nu} (A_{\varepsilon_0} +I)^{1/2} B^{-1} \Phi_2 =0$ at $J\times \Gamma$, since there holds
$(A_{\varepsilon_0} +I)^{1/2} B^{-1} = B^{-1} (A_{\varepsilon_0} +I)^{1/2} $ and
$\Phi_2 \in \zero{\sob{1/2}{p}{J;\lp{p}{\Omega}} \cap \lp{p}{J;\sob{2}{p}{\Omega}}}$, and
\begin{align*}
\del{}{\nu} A_{\varepsilon_0} ( c' - (A_{\varepsilon_0} +I)^{1/2} B^{-1} \Phi_2 )
& =
\del{}{\nu} ( A_{\varepsilon_0} c' - \Phi_2 ) + G^{3/4} \del{}{\nu} B^{-1} \Phi_2
\\
& \phantom{=\,}
+ \del{}{\nu} B^{-1} (A_{\varepsilon_0} +I)^{1/2} \Phi_2 = 0\quad\mbox{at}\;\,J\times \Gamma.
\end{align*}
We may thus estimate as follows.
\begin{align*}
\norm{ c' }{ \zero{Z_{2}}(J) } & \le C \norm{ c' }{ D(B) }\\
 & \le
C  \big( \norm{ c' - (A_{\varepsilon_0} +I)^{1/2} B^{-1} \Phi_2}{ D(B) } + \norm{ (A_{\varepsilon_0} +I)^{1/2} B^{-1} \Phi_2}{ D(B) } \big)
\\
& \le C \big( \norm{ B ( c' - (A_{\varepsilon_0} +I)^{1/2} B^{-1}  \Phi_2) }{ Y } + \norm{ (A_{\varepsilon_0} +I)^{1/2} \Phi_2 }{ Y } \big)
\\
& =  C \big( \norm{ L (G^{1/4} + (A_{\varepsilon_0}+I)^{1/2} )^{-1} F \big( c' - (A_{\varepsilon_0} +I)^{1/2} B^{-1}  \Phi_2\big) }{ Y }\\
& \quad+
\norm{ (A_{\varepsilon_0} +I)^{1/2} \Phi_2 }{ Y } \big)
\end{align*}
with $L:= B F^{-1} (G^{1/4} + (A_{\varepsilon_0}+I)^{1/2} )$. Note that, by reasoning as above, the operator
$H:=G^{1/4} + (A_{\varepsilon_0}+I)^{1/2}$ with domain
$D(H) := D(G^{1/4}) \cap D((A_{\varepsilon_0}+I)^{1/2}) = \sob{1/4}{2}{J;\lp{2}{\Omega}} \cap \lp{2}{J;\sob{1}{2}{\Omega}}$ is sectorial, invertible, and belongs to $\bip{Y}$
with power angle $\theta_H \le \pi/8$.

We now claim that
\begin{equation} \label{Lbound}
\norm{L \psi}{Y} \le \tilde{C} \norm{\psi}{Y},\quad \psi\in D(H).
\end{equation}
To see this, we rewrite $L$ as
\begin{align*}
L = B G^{1/4}  F^{-1} + B F^{-1} (A_{\varepsilon_0} +I)^{1/2}
=  B H F^{-1} + B F^{-1} \Big[(A_{\varepsilon_0} +I)^{1/2},F\Big] F^{-1}.
\end{align*}
Computation of the commutator results in
\begin{align*}
\Big[(A_{\varepsilon_0} +I)^{1/2},F\Big] = \Big[ (A_{\varepsilon_0} +I)^{1/2}, M \Big] G =: K_0 G
\end{align*}
where $K_0$ is a bounded operator. This can be seen by using the functional calculus for sectorial operators.
Thus $L$ can be represented as
\begin{align*}
L = B HF^{-1} + [ B F^{-1} ] K_0 [ G F^{-1}]
\end{align*}
which implies (\ref{Lbound}), in view of the regularity assumptions on
$\rho_0$, $\gamma_0$, and $\varepsilon_0$.

Combining the preceding estimates we obtain
\begin{equation*}
\norm{c'}{ \zero{Z_2(J)} } \le C_1\big( \norm{ H^{-1} F (c' - (A_{\varepsilon_0} +I)^{1/2} B^{-1} \Phi_2)}{ Y }
+ \norm{(A_{\varepsilon_0} +I)^{1/2} \Phi_2}{Y}\big).
\end{equation*}
Observe that $\tilde{c} := c' - (A_{\varepsilon_0} +I)^{1/2} B^{-1} \Phi_2$ solves
\begin{align*}
F \tilde{c} & = \Phi_1 - M G^{1/4} [ (A_{\varepsilon_0} +I)^{1/2} G^{1/2} B^{-1} ] G^{1/4} \Phi_2
+
G^{1/4} [A_{\varepsilon_0} G^{1/4} B^{-1}] G^{1/4} \Phi_2
\\
& \phantom{=\,} +
[ A_{\varepsilon_0}  (A_{\varepsilon_0}+I)^{1/2} B^{-1} ] \Phi_2  =: \Phi_3,
\end{align*}
where the operators inside the brackets $[ \ldots ]$ are bounded. Hence
the norm of $c'$ can be estimated by
\begin{equation} \label{eq:norm-c-1}
\norm{c'}{ \zero{Z_2(J)} } \le C_1\big( \norm{ H^{-1} \Phi_3 }{ Y } + \norm{(A_{\varepsilon_0} +I)^{1/2} \Phi_2}{Y}\big).
\end{equation}

{\bf c)} In order to estimate the first term on the right of \eqref{eq:norm-c-1} we will use duality relations.
In what follows $\du{\cdot}{\cdot}{}$ denotes the inner product
of $Y= \lp{2}{J;X} = \lp{2}{J;\lp{2}{\Omega}}$. Since $A_{\varepsilon_0} + I \in \hinfty{Y}$ is self-adjoint, we have
\begin{align*}
( (A_{\varepsilon_0} + I)^{1/2} )^* = ( A^*_{\varepsilon_0} + I)^{1/2} =  (A_{\varepsilon_0} + I)^{1/2},
\end{align*}
i.e. $(A_{\varepsilon_0} + I)^{1/2}$ is self-adjoint as well. One also readily verifies that $G^* = - G$ with domain
\begin{align*}
D(G^*) = {}^0\sob{1}{2}{[0,T];X}:=\{v \in \sob{1}{2}{[0,T];X}: v|_{t=T} = 0 \}.
\end{align*}
Using $(G^{1/4})^* = (G^*)^{1/4}$, cf. \cite[Proposition 5.1]{H1}, we are able to compute $D( (G^{1/4})^* )$
by complex interpolation leading to
$D( (G^{1/4})^*) = [Y, D(G^*)]_{1/4} = \sob{1/4}{2}{J;\lp{2}{\Omega}}$.
Observe that $(G^*)^{1/4}$ and $(A_{\varepsilon_0} + I)^{1/2}$ commute and belong to $\bip{Y}$
with power angles $\theta_{ (G^*)^{1/4}} \le \pi/8$ and $\theta_{ (A_{\varepsilon_0}+I)^{1/2} } =0$, respectively.
As above, we may conclude by the Dore-Venni theorem that $(G^*)^{1/4} + (A_{\varepsilon_0} + I)^{1/2}$ with natural domain
is invertible and belongs to $\bip{Y}$. In particular, there exists a constant $m > 0$
such that
\begin{align*}
\norm{(G^*)^{1/4} \psi}{ Y } + \norm{ (A_{\varepsilon_0} + I)^{1/2} \psi }{ Y } \le m
\norm{ (G^*)^{1/4} \psi + (A_{\varepsilon_0} +I)^{1/2} \psi }{ Y },
\end{align*}
for all $\psi \in D((G^*)^{1/4}) \cap D( (A_{\varepsilon_0}+I)^{1/2})$.
The operator $H=G^{1/4} + (A_{\varepsilon_0} + I)^{1/2}$, which was introduced in Step b), is densely defined, closed, and
invertible, which implies existence of $(H^{-1})^*$ and $(H^{-1})^* = (H^*)^{-1} \in \bdd{Y}{D(H^*)}$, where
$D(H^*) = \{ y' \in Y: \: \exists z \in Y: \: \du{H y}{y '}{} = \du{y}{z} \: \forall y \in D(H) \}$.
Furthermore, in view of the dense embedding $D(H) \hookrightarrow Y$, which entails
uniqueness of $z \in Y$, we obtain
\begin{align*}
D(H^*) = \{ y \in Y: \: (G^{1/4})^* y + (A_{\varepsilon_0} +I)^{1/2} y = z \in Y \}.
\end{align*}
The equation $(G^{1/4})^*y + (A_{\varepsilon_0} + I)^{1/2} y = z$ can be uniquely solved in
$Y$, that is for every $z \in Y$ there exists a unique  $y \in D( (G^*)^{1/4}) \cap D( (A_{\varepsilon_0} + I)^{1/2})$;
but this implies
\begin{align*}
D(H^*) = D( (G^*)^{1/4}) \cap D( (A_{\varepsilon_0} + I)^{1/2}) =
\sob{1/4}{2}{J;\lp{2}{\Omega}} \cap \lp{2}{J;\sob{1}{2}{\Omega}}.
\end{align*}

{\bf d)} We are now prepared to estimate the two terms on the right-hand side of \eqref{eq:norm-c-1}.
The first term can be rewritten as
\begin{equation} \label{eq:norm-c-2}
\begin{split}
\norm{ H^{-1} \Phi_3 }{ Y } & =\sup \{ |\du{ H^{-1} \Phi_3 }{ \psi' }|:  \|\psi'\|_{Y} \le 1 \}
\\
& = \sup \{ |\du{ \Phi_3 }{ \psi }| : \psi = (H^*)^{-1} \psi' \in D(H^*), \|H^* \psi \|_{Y} \le 1 \},
\end{split}
\end{equation}
so that it boils down to consider the inner product $\du{ \Phi_3}{ \psi }$ with $\psi \in D(H^*)$.
By definition of $\Phi_3$ we have
\begin{align*}
\du{ \Phi_3 }{  \psi} & = \du{ \Phi_1 }{ \psi } -
\du{ M [ (A_{\varepsilon_0} +I)^{1/2} G^{1/2} B^{-1} ] G^{1/4} \Phi_2 }{ (G^*)^{1/4} \psi }
\\
& \phantom{=\,}
+
\du{ [A_{\varepsilon_0} G^{1/4} B^{-1}] G^{1/4} \Phi_2 }{ (G^*)^{1/4} \psi }
+
\du{ [ A_{\varepsilon_0}  (A_{\varepsilon_0}+I)^{1/2} B^{-1} ] \Phi_2}{ \psi}
\end{align*}
and thus
\begin{align*}
|\du{ \Phi_3 }{  \psi}| \le |\du{ \Phi_1}{ \psi }| + ( C_1 \norm{ \Phi_2 }{ Y}
+
C_2 \norm{ G^{1/4} \Phi_2 }{Y} ) \norm{ (G^*)^{1/4} \psi }{ Y },
\end{align*}
where we used boundedness of the operators $M (A_{\varepsilon_0} +I)^{1/2} G^{1/2} B^{-1}$,
$A_{\varepsilon_0} G^{1/4} B^{-1}$, and $A_{\varepsilon_0}  (A_{\varepsilon_0}+I)^{1/2} B^{-1}$.

Next we deal with
the term $\du{ \Phi_1}{ \psi }$, where we aim at getting rid of one spatial derivative (it is actually the divergence operator
occurring in $A_{\gamma_0 - \gamma_1}$ etc.) by using the divergence theorem and the boundary conditions
$\del{}{\nu} \mu = \del{}{\nu} \mu_2 = 0$ and $\sk{\phi_2}{\nu} = 0$ (recall that $\sk{u_i}{\nu} = 0$, $i=1,2$),
so that all boundary integrals vanish. We obtain
\begin{equation} \label{est:Phi-1}
\begin{split}
\du{ \Phi_1 }{ \psi } & = \du{ \tfrac{\varepsilon_0}{\gamma_0}
\big ( G \phi_1 - \divv \phi_2 + A_{\gamma_0 - \gamma_1} \mu - A_{\gamma_1 - \gamma_2}\mu_2 \big)
+ \tfrac{\varepsilon_0^2}{\gamma_0} \nabla ( \tfrac{\gamma_0}{\varepsilon_0} ) \cdot
\nabla \mu }{ \psi}
\\
& = \du{ \tfrac{\varepsilon_0}{\gamma_0}  G^{3/4} \phi_1}{ (G^*)^{1/4}\psi}
+
\du{ \phi_2}{ \nabla ( \tfrac{\varepsilon_0}{\gamma_0} \psi )}
+
\du{ [\gamma_0 - \gamma_1] \nabla \mu}{ \nabla ( \tfrac{\varepsilon_0}{\gamma_0} \psi) }
\\
& \phantom{=\,}\,
-
\du{ [\gamma_1 - \gamma_2] \nabla\mu_2 }{ \nabla ( \tfrac{\varepsilon_0}{\gamma_0}\psi)}
+
\du{ ( \tfrac{\varepsilon_0^2}{\gamma_0}  \nabla ( \tfrac{\gamma_0}{\varepsilon_0} ) \cdot \nabla \mu}{ \psi) }
\end{split}
\end{equation}
and thus
\begin{align*}
|\du{ \Phi_1 }{ \psi }|
& \le C_1 \big( \norm{ G^{3/4} \phi_1 }{ Y } + \norm{\phi_2}{Y}
+
\norm{ \gamma_0 - \gamma_1 }{ \cont{}{}{J \times \ol{\Omega}} }
\norm{ \nabla \mu }{ Y }
\nonumber
\\
& \phantom{\le \,}
\,+
\norm{ \gamma_1- \gamma_2 }{ \lp{\infty}{J;\lp{2}{\Omega}} }
\norm{ \nabla \mu_2 }{ \lp{2}{J;\lp{\infty}{\Omega;\R^n}} }
+
C_2 \norm{ \nabla \mu }{ Y } \norm{ \psi }{ Y }\nonumber
\\
& \le \ldots + C_2 T^{1/4} \norm{ \nabla \mu }{ Y } \norm{ \psi }{ D( (G^*)^{1/4}) },
\end{align*}
where H\"older's inequality and the embedding $\psi \in D(G^*) \hookrightarrow \lp{4}{J;\lp{2}{\Omega}}$ entered.

Taking the supremum w.r.t. $\norm{ H^* \psi }{ Y } \le 1$, we get
\begin{align}
\norm{c'}{ \zero{Z_2(J)} } & \le C \Big(
\norm{ G^{3/4} \phi_1 }{ Y } + \norm{\phi_2}{Y}
+
\norm{ \gamma_0 - \gamma_1 }{\infty} \norm{ \nabla \mu }{ Y } \nonumber\\
& \quad +
\norm{ \gamma_1- \gamma_2 }{ \lp{\infty}{J;\lp{2}{\Omega}} }
\norm{ \nabla \mu_2 }{ \lp{2}{J;\lp{\infty}{\Omega;\R^n}} }
+
T^{1/4} \norm{ \nabla \mu }{ Y }\nonumber\\
& \quad +
\norm{ \Phi_2 }{ Y} + \norm{ G^{1/4} \Phi_2 }{Y} \Big)
+
\norm{ (A_{\varepsilon_0} + I)^{1/2} \Phi_2 }{ Y }\nonumber
\\
 & \le\,
C \big(
\norm{ G^{3/4} \phi_1 }{ Y } + \norm{\phi_2}{Y}
+
( \norm{ \gamma_0 - \gamma_1 }{\infty} + T^{1/4} ) \norm{\mu}{ Z_\mu(J) }
\nonumber\\
&\quad +
\norm{ \gamma_1- \gamma_2 }{ \lp{\infty}{J;\lp{2}{\Omega}} }
+
\norm{ \Phi_2 }{ Z_\mu(J) } \big). \label{est:c-1}
\end{align}
Before continuing with this estimate, we collect some embeddings
which are frequently used below. Moreover, as compositions of nonlinear functions with
$c_i$, $u_i$, $\rho_i$ occur, we will also need some preliminary considerations to this matter.

{\bf e)} We first remind the reader that the spaces ${\mathcal{W}_i(J)}$ were defined in \eqref{embedd:2}.
Set $p_1 := 2p/(p-2)$, $\alpha := 1/4(1-n/p)$, and let $\delta \in
[0,1/4)$ and $\theta \in (0,1)$.
Since $p>\hat{p}=\max\{4,n\}$, we have $\alpha > 0$,
$c_0 \in \sobb{4-4/p}{p}{\Omega} \hookrightarrow \cont{2}{}{\ol{\Omega}}$,
$\rho_0 \in \sob{3}{p}{\Omega} \hookrightarrow \cont{2}{}{\ol{\Omega}}$, and the embeddings
\begin{equation} \label{embedd:3}
\begin{split}
& c_0 - c_1 \in \zero{\mathcal{W}_2(J)}, \quad \rho_0 - \rho_1 \in
\zero{\mathcal{W}_3(J)}, \quad \sob{1}{2}{\Omega} \hookrightarrow \lp{p_1}{\Omega},
\\
&u \in \zero{Z_1(J)} \hookrightarrow \zero{\cont{\delta + 1/2}{}{J;\lp{2}{{\Omega;\R^n}}}} \cap
\zero{\cont{\delta+1/4}{}{J;\lp{p_1}{\Omega;\R^n}}} \cap \sob{1/2}{2}{J;\sob{2}{2}{\Omega;\R^n}},
\\
& c \in \zero{Z_2(J)} \hookrightarrow
\sob{\frac{3}{4} \theta}{2}{J;\sob{3(1-\theta)}{2}{\Omega}}\\
& \quad\quad\quad\hookrightarrow \zero{\sob{1/2}{2}{J;\sob{1}{2}{\Omega}}} \cap
\zero{\sob{1/4+\alpha}{2}{J;\sob{1}{p_1}{\Omega}}} \cap \zero{\cont{\delta}{}{J;\lp{2}{{\Omega}}}},
\\
& \mu \in \zero{Z_\mu(J)} \hookrightarrow \lp{4}{J;\lp{2}{\Omega}} \cap \lp{2}{J;\lp{p_1}{\Omega}},
\\
& \rho \in \zero{Z_3(J)} = \zero{\sob{2}{2}{J;\lp{2}{\Omega}}} \cap \zero{\cont{1}{}{J;\sob{1}{2}{\Omega}}}
\cap \zero{\cont{}{}{J;\sob{2}{2}{\Omega}}}. 
\end{split}
\end{equation}
In what follows we will extensively use the inequality (here $X$ belongs to the class $\mathcal{HT}$)
\begin{align} \label{eq:sob-T}
\norm{\phi}{ \zero{\sob{\gamma}{p}{J;X}} }
\le
T^{1-\gamma} \norm{ \phi }{ \zero{\sob{1}{p}{J;X}} }, \quad
\forall \phi \in \sob{1}{p}{J;X}, \quad p \in (1,\infty), \quad \gamma \in [0,1).
\end{align}
To treat differences of nonlinear functions, we need some obvious, but very useful estimates. Let $\alpha\in \cont{2}{}{\R^{n^2+n+2};\R^m}$, $m\in \N$,
$(u_i,c_i,\mu_i) \in \Sigma_T$, and $\rho_i = L[u_i] \rho_0$, $i=1,2$.
Then there exist constants $k,C > 0$ which are independent of the triples $(u_i,c_i,\mu_i)$, $i=1,2$, and $T$ such that
\begin{multline} \label{diff:nl}
\sum_{l=0}^2 \left| \nabla^l \big[\alpha(u_i,c_i,\rho_i,\nabla( u_i, c_i, \rho_i)) (t,x)
- \alpha(u_j,c_j,\rho_j,\nabla( u_j, c_j, \rho_j))(t,x)\big]\right |
\\
\le k_0  \sum_{l=0}^3 |\nabla^l (u_i,c_i,\rho_i)(t,x) - \nabla^l (u_j,c_j,\rho_j)(t,x) |
\end{multline}
with $k_0 := C \max\{ \sup_{ |y| \le k } | \alpha'(y) |,
\sup_{|y| \le k} |\alpha''(y)|, \sup_{|y| \le k} |\alpha'''(y)| \}$ and $i,j \in \{1,2,0\}$. This is possible due to
the embedding $\mathcal{Z}(J) \hookrightarrow \cont{}{}{J;\cont{2}{}{\ol{\Omega};\R^{n+2}}}$
and the estimates \eqref{est:uniform}, \eqref{est:rho}.

After these preparations we are now able to estimate differences of the nonlinear functions in \eqref{est:c-1}
in the function spaces under consideration. We have (with some $\delta\in (0,1/4)$)
\begin{align*}
\norm{ \gamma_0 - \gamma_1 }{ \zero{\cont{}{}{J \times \ol{\Omega}}} }
& \le k_0 \big( \norm{ c_0 - c_1 }{ \zero{\cont{}{}{J \times \ol{\Omega}}} }
+
\norm{ \rho_0 - \rho_1 }{ \zero{\cont{}{}{J \times \ol{\Omega}}} } \big)
\le C ( T^{1/2} + T ),
\\
\norm{ \gamma_1 - \gamma_2 }{ \zero{\cont{}{}{J;\lp{2}{\Omega}}} }
& \le k_0 ( \norm{ c }{ \zero{\cont{}{}{J;\lp{2}{\Omega}}} }
+
\norm{ \rho }{ \zero{\cont{}{}{J;\lp{2}{\Omega}}} } )
\le
\tilde{k}_0 ( T^{\delta} \norm{ c }{ \zero{Z_2(J)}}
+
T \norm{ \rho }{ \zero{ Z_3(J) } } ).
\end{align*}
In a similar manner we deal with $G^{3/4} \phi_1$ and $\phi_2$ occurring in \eqref{est:c-1},
\begin{align*}
& \norm{G^{3/4} \phi_1}{Y} \le C \norm{\phi_1}{ \zero{\sob{3/4}{2}{J;\lp{2}{\Omega}}} }
\le C (
T^{1-\beta} \norm{\rho_0 - \rho_1}{ \zero{\cont{1}{}{J;\cont{}{}{\ol{\Omega}}}} } \norm{c}{ \zero{Z_2(J)} }
\\
& \quad +
T^{1-\beta} \norm{ \rho }{ \zero{\cont{1}{}{J;\lp{2}{\Omega}}} }
\norm{c_2}{ \sob{3/4}{2}{J;\cont{}{}{\ol{\Omega}}} } )
\le
C T^{1-\beta} ( \norm{c}{ \zero{Z_2(J)} } + \norm{ \rho }{ \zero{ Z_3(J) } } ), \: \tfrac{3}{4} < \beta <1,
\\
& \norm{ \phi_2 }{ Y } \le k_0 T^{1/2}  \big(
T^{1/2+\delta} \norm{  u }{ \zero{ Z_1(J) } } + T^{\delta} \norm{  c }{ \zero{ Z_2(J) } } + T \norm{  \rho }{ \zero{ Z_3(J) } }
\big).
\end{align*}
Looking at the definition of  $\Phi_2$, one realizes that it consists of the highest
order term $-A_{\varepsilon_1 - \varepsilon_0} c$ and lower order terms that provide a factor
of the type $T^{a}$, $a > 0$, due to more time regularity.
Although $\rho$ is measured w.r.t.\ space in $\sob{2}{2}{\Omega}$,
the term $\rho_1 \varepsilon_1 \nabla c_1 \cdot \nabla \rho$ is still of lower order,
since we have continuity w.r.t. the time variable, cf. the definition of $Z_3(J)$ and estimate
\eqref{est:diff-rho}. Using \eqref{diff:nl} and the embeddings \eqref{embedd:2}, \eqref{embedd:3}
one finds after some tedious computations the following estimate
\begin{multline*}
\norm{ \Phi_2 }{ Z_\mu(J) } \le
C_1 \norm{ \varepsilon_0 - \varepsilon_1 }{ \cont{1/4}{}{J;\cont{}{}{\ol{\Omega}}} \cap \cont{}{}{J;\cont{1}{}{\ol{\Omega}} } }
\norm{c}{ \zero{Z_2(J)} }
+
C_2 \norm{  \nabla( \varepsilon_1 - \varepsilon_2) }{ Z_\mu(J) } \cdot
\\
\norm{ \nabla c_2 }{ \cont{1/4+\beta}{}{J;\cont{}{}{\ol{\Omega}}} \cap \cont{}{}{J;\cont{1}{}{\ol{\Omega}}} }
+
C_2 \norm{  \varepsilon_1 - \varepsilon_2 }{ \lp{2p/(p-2)}{J;\sob{1}{2}{\Omega}} }
\norm{ \Delta c_2 }{ \lp{p}{J;\cont{1}{}{\ol{\Omega}}} }
\\
+
C_2 \norm{  \varepsilon_1 - \varepsilon_2 }{ \sob{1/4}{2}{J;\lp{2p/(p-2)}{\Omega}} }
\norm{ \Delta c_2 }{ \cont{1/4 + \beta}{}{J;\lp{p}{\Omega}} }
+
C_3( \norm{\rho}{ Z_\mu(J) } + \norm{ \nabla \rho }{ Z_\mu(J) }
\\
+ \norm{ c }{ Z_\mu(J) } + \norm{ \nabla c }{ Z_\mu(J) } )
\le C_1 k_0 (T^{1/4+\beta} + T^{3/4} + T^{1/2}) \norm{c}{ \zero{Z_2(J)} }
\\
+
C_2 k_0 (T^{1/4} \norm{c}{ \zero{Z_2(J)} }+ (T^{1/2} + T^{1/4}) \norm{\rho}{ \zero{Z_3(J)} } )
\\
+ C_3 (T^{1/4} \norm{c}{ \zero{Z_2(J)} } + (T^{1/2} + T^{1/4}) \norm{\rho}{ \zero{Z_3(J)} } ).
\end{multline*}
By means of \eqref{est:diff-rho} we may then replace all difference terms involving
$\rho = \rho_1 - \rho_2$ by corresponding terms with $u=u_1- u_2$, thereby obtaining
\begin{align*}
\norm{c'}{\zero{Z_2(J)} } \le \kappa'_2(T) ( \norm{u}{ \zero{Z_1(J)} }  + \norm{c}{ \zero{Z_2(J)} } + \norm{\mu}{ Z_\mu(J) } ).
\end{align*}
Finally, as already mentioned, we easily find an estimate for $\mu'$ in $Z_\mu(J)$ by using the equation $\mu' =A_{\varepsilon_0} c' - \Phi_2$; this gives
\begin{align*}
\norm{\mu'}{ Z_\mu(J) } & \le \norm{ A_{\varepsilon_0} c' }{ Z_\mu(J) } + \norm{\Phi_2}{ Z_\mu(J)}
\le C \norm{ c ' }{ \zero{Z_2(J)} } + \norm{\Phi_2}{ Z_\mu(J)}
\\
& \le C \kappa'_2(T) ( \norm{u}{ \zero{Z_1(J)} } + \norm{ c }{ \zero{Z_2(J)} } + \norm{ \mu }{ Z_\mu } ),
\end{align*}
which completes the proof by summarizing both estimates.
\eproof

\section{The final step in the proof of Theorem \ref{theo:main:1}} \label{sec:PrThm}
We recall that in Section \ref{sec:3} the original problem
\eqref{eq:nsch-1}-\eqref{bc:3} was equivalently rewritten as
\begin{align*}
\mathcal{ L}w =\mathcal{F}(w,\rho), \quad \rho = L[u]\rho_0,
\end{align*}
where $w=(u,c,\mu)$, $\rho=L[u]\rho_0$ is the solution of the continuity equation, and where the linear operator $\mathcal{ L}$
and the nonlinear operator $F(w,\rho)$ are as defined as in Section \ref{sec:3}.
 By Lemma \ref{lem:rho}, we know that for any $T\in (0,T_0]$, $J=[0,T]$, and $T_0>0$ being fixed, we have
\begin{align*}
\norm{\rho}{ \mathcal{Z}_3(J)}\le \varsigma \quad \mbox{for all} \;u \in \mathcal{Z}_{1,\Gamma}(J),
\, \rho_0 \in \sob{3}{p}{\Omega}, \, \rho_0>0,
\end{align*}
where we set $\mathcal{Z}_{1,\Gamma}(J) := \{u \in \mathcal{Z}_1(J): \: \ska{u}{\nu}{\Gamma} \ge 0 \; \mbox{and}\;
\norm{u - \ol{u}}{ \zero{\mathcal{Z}_1(J)} } \le 1$\}, $\ol{u}\in \mathcal{Z}_1(J_0)$ is a fixed function, and $\varsigma>0$ 
does not depend on $T$ and $u$. Furthermore, due to the Theorems \ref{theo:ch} and \ref{theo:stokes:3} we have maximal regularity for
the associated linear problem, that is the operator $\mathcal{L}$ is an isomorphism from the desired
class of maximal regularity to the corresponding space of data. More precisely, we have
\begin{align*}
\mathcal{L} & \in \liss{\mathcal{Z}_{1,\mathcal{B}}(J) \times \mathcal{Z}_2(J)\times \mathcal{Z}_\mu(J)}{
\mathcal{D}_1(J) \times \mathcal{D}_2(J) \times \mathcal{D}_\mu(J)},
\\
\mathcal{D}_1(J) & := \{ \varphi \in \mathcal{X}_{1,n,\Gamma}(J) \times
\mathcal{Y}_{0,d}(J;\R^n) \times \mathcal{Y}_{0,s}(J)
\times \mathcal{Y}_{1,s}(J;\R^n) \times \sobb{4-2/p}{p}{\Omega;\R^n}: \\
& \qquad \varphi \text{ fulfils 4.-7. of Theorem \ref{theo:stokes:3}} \},
\\
\mathcal{D}_2(J) & := \{ \varphi \in \mathcal{X}_2(J) \times \mathcal{Y}_1(J) \times \sobb{4-4/p}{p}{\Omega}:
\varphi \text{ fulfils 4. of Theorem \ref{theo:ch} } \},
\\
\mathcal{D}_\mu(J) & := \{ \varphi \in \mathcal{X}_\mu(J) \times \mathcal{Y}_\mu(J):
\varphi \text{ fulfils 4. of Theorem \ref{theo:ch} } \}.
\end{align*}
Using this property and Lemma \ref{lem:rho} as well as the assumptions in Theorem
\ref{theo:main:1}, it is not difficult to verify that $\mathcal{F}$ (see \eqref{def:F} for its definition), maps
$\mathcal{Z}_{1,\mathcal{B}}(J) \times \mathcal{Z}_2(J)\times \mathcal{Z}_\mu(J)\times \mathcal{Z}_3(J)$ to
$\mathcal{D}_1(J)\times \mathcal{D}_2(J) \times \mathcal{D}_\mu(J)$, see also the estimates below. Now, for any $w$ from the set
\begin{align*}
\begin{split}
\Sigma_T = \{ &(u,c,\mu) \in \mathcal{Z}_1(J) \times
\mathcal{Z}_2(J)\times \mathcal{Z}_\mu(J):\;\,(u,\partial_t u,c,\mu)(0) = (u_0,u_\bullet,c_0,\mu_0),\\
& \;u=0\;\mbox{on}\;J\times \Gamma_d,\;\ska{ u }{ \nu }{} = 0\;\mbox{and}\;
\mathcal{Q} \mathcal{D}(u) \cdot \nu=0\;\mbox{on}\;J\times\Gamma_s,
\\
& \;\partial_\nu c=\partial_\nu \mu=0\;\mbox{on}\;J\times \Gamma,\;\mbox{and}\\
& \;
\norm{ (u,c,\mu) - (\ol{u},\ol{c},\ol{\mu}) }{
{\mathcal{Z}_1}(J) \times {\mathcal{ Z}_2}(J)\times \mathcal{Z}_\mu(J) } \le 1 \}
\end{split}
\end{align*}
(see below for the definition of $(\ol{u},\ol{c},\ol{\mu})$),
we have $w\in \mathcal{Z}_{1,\mathcal{B}}(J) \times \mathcal{Z}_2(J)\times \mathcal{Z}_\mu(J)$, by definition of
$\Sigma_T$. Thus the problem $\mathcal{L}w' =\mathcal{F}(w,\rho)$ with $\rho=L[u]\rho_0$
has a unique solution $w'\in \mathcal{Z}_{1,\mathcal{B}}(J) \times \mathcal{Z}_2(J)\times \mathcal{Z}_\mu(J)$, that is
the fixed point mapping $\mathcal{G}$, which was introduced in (\ref{eq:fpe-mapping}), is
well-defined as a mapping from $\Sigma_T$ to $\mathcal{Z}_{1,\mathcal{B}}(J) \times \mathcal{Z}_2(J)\times \mathcal{Z}_\mu(J)$. We will
show that for sufficiently small $T$ the fixed point mapping $\mathcal{G}$ leaves $\Sigma_T$
invariant and becomes a strict contraction in $Z_1(J)\times Z_2(J)\times Z_\mu(J)$.

We now choose the reference function $\ol{w}:=(\ol{u},\ol{c},\ol{\mu})$ as the unique solution of
\begin{align} \label{def:centre}
\mathcal{L} \ol{w} = (\mathcal{F}_1(w_0,\rho_0),u_0,\mathcal{F}_2(w_0,\rho_0),c_0,\mathcal{F}_\mu(w_0,\rho_0) )
\end{align}
on the time interval $J_0=[0,T_0]$ with $w_0 := (u_0,c_0,\mu_0)$ and $\mu_0$ as in (\ref{mu0def}).
Note that $\mathcal{F}_i(w_0,\rho_0)$, $i=1,2,\mu$,
comprises only terms of lower order so that the right-hand side
of (\ref{def:centre}) belongs to $\mathcal{D}_1(J_0) \times \mathcal{D}_2(J_0) \times \mathcal{D}_\mu(J_0)$, in particular
all compatibility conditions are satisfied. Theorem \ref{theo:ch} and \ref{theo:stokes:3}
guarantee existence and uniqueness of $(\ol{u},\ol{c},\ol{\mu})$ in
$\mathcal{Z}_{1,\mathcal{B}}(J_0) \times \mathcal{Z}_2(J_0) \times \mathcal{Z}_\mu(J_0)$
for any $T_0\in (0,\infty)$.

We recall that for any function space $Y(J)$ with the properties $\mathcal{Z}_1(J) \times \mathcal{Z}_2(J)
\times \mathcal{Z}_\mu(J)\hookrightarrow Y(J)$, $\norm{{v}}{Y(J)}\le \norm{{v}}{Y(J_0)}$ for all $v\in Y(J_0)$, 
and such that there exists a bounded extension operator $E_+$ from $\zero{ Y(J) }$ to
$\zero{ Y }(J_0)$ with a bound $\|E_+\|$ which is uniform w.r.t.\ $T\in (0,T_0]$, we have
$\norm{(u,c,\mu)}{ Y(J) } \le k$ for any $(u,c,\mu) \in \Sigma_T$ where $k$ does not depend on $T$, see
the remarks after Lemma \ref{lem:rho}. Moreover, as in the proof of Lemma \ref{lem:diff-ch},
all subsequent constants $C$, $C_i$, $C_{ij}$, $i,j \in \N$,
which may differ from line to line, are always independent of the unknowns and $T\in (0,T_0]$,
which is again due to working in function spaces with vanishing temporal traces at $t=0$.

{\bf Step 1: Contraction.} Let $w_i := (u_i,c_i,\mu_i) \in \Sigma_{T}$, $i=1,2$,
$\rho_i = L[u_i]\rho_0$. In this part of the proof, we denote by $w'=(u',c',\mu')$ the unique solution of
\begin{align*}
\mathcal{L} w' =
(\mathcal{F}_1(w_1,\rho_1) - \mathcal{F}_1(w_2,\rho_2),0,\mathcal{F}_2(w_1,\rho_1) - \mathcal{F}_2(w_2,\rho_2),0,
\mathcal{F}_\mu(w_1,\rho_1) - \mathcal{F}_\mu(w_2,\rho_2)).
\end{align*}
As in Section \ref{sec:WCH}, we set again $u:=u_1 - u_2$, $c:=c_1 - c_2$, $\mu:=\mu_1 - \mu_2$, and
$\rho := \rho_1 - \rho_2$. To obtain an estimate for $w'$ in $\zero{Z_1(J)} \times \zero{Z_2(J)} \times {}_0Z_\mu(J)$,
we just apply the maximal regularity result Theorem \ref{theo:stokes:2} to the
equation for $u'$ and Lemma \ref{lem:diff-ch} to the problem for $(c',\mu')$ leading to
\begin{align*}
\norm{ c' }{ \zero{Z_2(J)} } + \norm{ \mu' }{\zero{ Z_\mu(J) }}
& \le \kappa_2(T) \norm{ (u,c,\mu) }{ \zero{ Z_1(J)} \times \zero{ Z_2(J)} \times \zero{ Z_\mu(J)}}
\end{align*}
with $\kappa_2(T) \to 0$ as $T \to 0$, and
\begin{align}
\norm{u'}{ \zero{Z_{1}(J)} }
 & \le M_1 \left \{
\norm{ \mathcal{B}(D) c') }{ X_1(J) }
+
\norm{ \mathcal{C} \mu' }{ X_1(J) }
+
\norm{ B_{low}(w_1,\rho_1) - B_{low}(w_2,\rho_2)}{ X_1(J) }
\right.
\nonumber\\
&  +\norm{ B_1(w_1,\rho_1)u_1 - B_1(w_2,\rho_2)u_2}{ X_1(J) }
+
\norm{ B_2(w_1,\rho_1)c_1 - B_2(w_2,\rho_2)c_2 }{ X_1(J) }
\nonumber\\
& \left.
+\norm{ B_3(w_1,\rho_1)\rho_1 - B_3(w_2,\rho_2)\rho_2}{ X_1(J) }
+
\norm{ B_{\mu}(w_1,\rho_1)\mu_1 - B_{\mu}(w_2,\rho_2)\mu_2}{ X_1(J) }
\right \}
\nonumber\\
& \le M_1 C \kappa_2(T) \norm{ (u,c,\mu) }{ \zero{ Z_1(J)} \times \zero{ Z_2(J)} \times \zero{Z_\mu(J) }}
\nonumber\\
& \quad+
M_1 \{ \norm{ B_{low}(w_1,\rho_1) - B_{low}(w_2,\rho_2)}{ X_1(J) } + \ldots \}, \label{claimcont}
\end{align}
where we set $X_1(J) := \sob{1/4}{2}{J;\lp{2}{\Omega;\R^n}} \cap
\lp{2}{J;\sob{1}{2}{\Omega;\R^n}}$, which coincides with the $n$-dimen\-sional case of $Z_\mu(J)$.
Moreover, in view of Lemma \ref{lem:diff-rho}, each difference $\rho=\rho_1 - \rho_2$ occurring above
can be estimated by means of $u=u_1 - u_2$,
\begin{align*}
\norm{\rho}{ \zero{Z_3(J)} } \le \kappa_1(T) \norm{u}{ \zero{Z_1(J)} }.
\end{align*}
Subsequently, it is crucial that $M_1$, the operator norm of ${\cal L}_1^{-1}$,
is independent of the time interval $J = [0,T]$, $T\in (0,T_0]$, but might depend on
the fixed time $T_0$. This is due to the fact that we have spaces with zero initial data, which comes from
looking at differences of functions with the same initial data. This property will also be used in the estimates
below, where constants occur due to embedding and interpolation inequalities.

We now continue to estimate the second term on the right of (\ref{claimcont}).
Note that $B_i(w_0,\rho_0) \equiv 0$, $i=1,2,\mu$ (see \eqref{def:B} and the remarks below),
and recall the embeddings \eqref{embedd:1},
\eqref{embedd:2}, and \eqref{embedd:3}. For $B_\mu$, one easily finds
\begin{multline*}
\norm{ B_{\mu}(w_1,\rho_1)\mu_1 - B_{\mu}(w_2,\rho_2)\mu_2}{ X_1(J) } \le
\norm{ \rho_0 \nabla c_0 - \rho_1 \nabla c_1 }{ \cont{1/4+\beta}{}{J;\cont{}{}{\ol{\Omega}}}}
\norm{ \mu }{ \sob{1/4}{2}{J;\lp{2}{\Omega}}}
+
\\
\norm{ \rho_0 \nabla c_0 - \rho_1 \nabla c_1 }{ \cont{}{}{J;\cont{}{}{\ol{\Omega}}}}
\norm{ \mu }{ \lp{2}{J;\sob{1}{2}{\Omega}}}
+
\norm{ \rho_0 \nabla c_0 - \rho_1 \nabla c_1 }{ \cont{}{}{J;\cont{1}{}{\ol{\Omega}}}}
\norm{ \mu }{ \lp{2}{J;\lp{2}{\Omega}}}
+
\\
\norm{ \rho_1 \nabla c_1 - \rho_2 \nabla c_2 }{ \lp{p_1}{J;\lp{2}{\Omega;\R^n}} }
\norm{ \mu_2 }{ \lp{p}{J;\cont{1}{}{\ol{\Omega}}} }
+
\norm{ \rho_1 \nabla c_1 - \rho_2 \nabla c_2 }{ \lp{2}{J;\sob{1}{2}{\Omega;\R^n}} }
\norm{ \mu_2 }{ \cont{}{}{J \times \ol{\Omega}} }
\\
+ \norm{ \rho_1 \nabla c_1 - \rho_2 \nabla c_2 }{ \sob{1/4}{2}{J;\lp{p_1}{\Omega;\R^n}} }
\norm{ \mu_2 }{ \cont{1/4+\beta}{}{J;\lp{p}{\Omega}} } ,
\end{multline*}
where again $p_1 := 2p/(p-2)$, and this can be further estimated from above by
\begin{multline*}
\big( C_1 ( T^{1/4-\beta} + T^{1/2} ) + C_2 (T^{1/2} + T) + C_3 T^{1/4} \big) \norm{ \mu }{ Z_\mu(J) }
+
C_4 T^{1/p_1 } \norm{c}{ \zero{Z_2(J)} }
+
\\
C_4 T^{1+1/p_1} \norm{\rho}{ \zero{Z_3(J)} }
+
C_5 ( T^{1/4} \norm{c}{ \zero{Z_2(J)} } + T^{1/2+1} \norm{\rho}{ \zero{Z_3(J)} })
+
C_6 ( T^{\alpha} \norm{c}{ \zero{Z_2(J)} } + T  \norm{\rho}{ \zero{Z_3(J)} } ).
\end{multline*}
Here we get uniform estimates (w.r.t.\ elements of $\Sigma_T$) for the norm of terms of the 
structure 'function - initial value' by using the reference functions from the definition of $\Sigma_T$, 
e.g. (with $\rho(\ol{u}):=L[\ol{u}]\rho_0$)
\begin{align*}
\norm{ \rho_1 \nabla c_1 - \rho_0 \nabla c_0 }{ \cont{1/4+\beta}{}{J;\cont{}{}{\ol{\Omega}}}}\,\le & 
\,\norm{ \rho_1 \nabla c_1 - \rho(\ol{u}) \nabla \ol{c} }{ \cont{1/4+\beta}{}{J;\cont{}{}{\ol{\Omega}}}}
\\ 
& + \norm{ \rho (\bar{u}) \nabla \bar{c} - \rho_0 \nabla c_0 }{ \cont{1/4+\beta}{}{J;\cont{}{}{\ol{\Omega}}}};
\end{align*}
the second term is independent of $w$, and for the first one we use the condition $\norm{ (u,c,\mu) - (\ol{u},\ol{c},\ol{\mu}) }{
{\mathcal{Z}_1}(J) \times {\mathcal{ Z}_2}(J)\times \mathcal{Z}_\mu(J) } \le 1$ from the definition of $\Sigma_T$.

Turning to the $B_{low}$-term we have
\begin{align*}
\exists (\tilde{u},\tilde{c},\tilde{\mu}) \in \Sigma_T\;\mbox{s.t.} \;
B_{low}(w_1,\rho_1) - B_{low}(w_2,\rho_2)
=
\tilde{B}_{low} ( \tilde{u},\tilde{c},\tilde{\rho} ) (u, c, \rho, \nabla u, \nabla c)\;\in X_1(J)
\end{align*}
where $\tilde{\rho} = L[\tilde{u}] \rho_0$ and
$\tilde{B}_{low} \in \bdd{ \mathcal{Z}(J) }{  \mathcal{M}_\mu(J;\R^{n+2} \times \R^{n \times (n+1)}) }$
with $\mathcal{M}_\mu$ denoting the multiplier space
$\mathcal{M}_\mu(J;E) = \cont{1/4 + \beta}{}{J;\cont{}{}{\ol{\Omega};E}} \cap \cont{}{}{J;\cont{1}{}{\ol{\Omega};E}}$
of $Z_\mu(J)$. Note that $\sob{1/4}{2}{J;\lp{2}{\Omega}}=\sobb{1/4}{2}{J;\lp{2}{\Omega}}$, 
so that here we can use well-known results on pointwise multipliers for vector-valued Besov spaces,
see e.g.\ \cite{amannM}. We obtain
\begin{multline*}
\norm{ B_{low}(w_1,\rho_1) - B_{low}(w_2,\rho_2)}{ X_1(J) }
\le C_1
\norm{ u }{ \sob{1/4}{2}{J;\sob{1}{2}{\Omega;\R^{n}}} \cap \lp{2}{J;\sob{2}{2}{\Omega;\R^{n}}} }
\\
+C_2 \norm{ c }{ \sob{1/4}{2}{J;\sob{1}{2}{\Omega}} \cap \lp{2}{J;\sob{2}{2}{\Omega}} }
+
C_3 \norm{ \rho }{ \sob{1/4}{2}{J;\lp{2}{\Omega}} \cap \lp{2}{J;\sob{1}{2}{\Omega}} }
\\
\le C_1 (T^{3/4} + T^{1/2} ) \norm{u}{ \zero{Z_1(J)} }
+
C_2 T^{1/4} \norm{ c }{ \zero{Z_2(J)} } + C_3 ( T^{1/2} + T^{3/4} T^{1/2} ) \norm{\rho}{ \zero{Z_3(J)} }.
\end{multline*}
Next, we write the difference $B_1(w_1,\rho_1)u_1 - B_1(w_2,\rho_2)u_2$ as
\begin{equation*}
[B_1(w_1,\rho_1) - B_1(w_0,\rho_0)]u+[B_1(w_1,\rho_1) - B_1(w_2,\rho_2)]u_2
\end{equation*}
and get
\begin{multline*}  
\norm{[B_1(w_1,\rho_1) - B_1(w_0,\rho_0)]u}{ X_1(J) } \le
C_1 \norm{ \rho_1 - \rho_0 }{ \mathcal{M}(J;\R) } \norm{\del{}{t} u}{ Z_\mu(J) }
\\
+C_2 \norm{(\rho_1 - \rho_0,c_1 - c_0)}{ \mathcal{M}(J;\R^2) } \norm{\nabla^2 u}{ Z_\mu(J) }
+
C_3 \norm{ \nabla (\rho_1-\rho_0,c_1-c_0) }{ \mathcal{M}(J;\R^{n \times 2}) } \norm{ \nabla u }{ Z_\mu(J) }
\\
\le C_1 (T^{3/4-\beta}  + T^{1/2} ) \norm{u}{ \zero{Z_1(J)} }
+
C_2 ( T^{1/4} + T^{1/4+\beta} ) \norm{u}{ \zero{Z_1(J)} }
+
C_2 ( T^{3/4} + T^{1/2} ) \norm{u}{ \zero{Z_1(J)} }
\end{multline*}
and
\begin{multline*}  
\norm{[B_1(w_1,\rho_1) - B_1(w_2,\rho_2)]u_2}{ X_1(J) } \le
C_1 ( \norm{ \rho }{ \sob{1/4}{2}{\lp{2}{\Omega}} }
\norm{ \del{}{t} u_2 }{ \cont{1/4+\beta}{}{J;\cont{}{}{\ol{\Omega};\R^n}} }
\\
+\norm{ \rho }{ \cont{}{}{J;\sob{1}{2}{\Omega}} } \norm{ \del{}{t} u_2 }{ \lp{2}{J;\cont{1}{}{\ol{\Omega};\R^n}}}
)
+
C_2 ( \norm{(\rho,c)}{ \lp{p_1}{J;\sob{1}{2}{\Omega}} } \norm{u_2}{ \lp{p}{J;\cont{3}{}{\ol{\Omega};\R^n}} }
\\
+
\norm{(\rho,c)}{ \sob{1/4}{2}{J;\lp{2}{\Omega}} }
\norm{u_2}{ \cont{1/4+\beta}{}{J;\cont{2}{}{\ol{\Omega};\R^n}}} )
+
C_3 \norm{ \nabla (\rho,c) }{ Z_\mu(J) } \norm{ \nabla u_2 }{ \mathcal{M}(J;\R^{n \times n}) } )
\\
\le C_1 (T^{1+3/4}  + T ) \norm{\rho}{ \zero{Z_3(J)} } + C_2 ( T^{1/p_1} + T^{1/2} ) \norm{c}{ \zero{Z_2(J)} }
\\
+C_2 (T^{1/p_1 + 1} + T^{3/4 + 1/2} ) \norm{\rho}{\zero{Z_3(J)}}
+
C_3 T^{1/4} \norm{c}{ \zero{Z_2(J)} }
+
C_3 ( T^{1/2} + T^{1/2+3/4} ) \norm{\rho}{ \zero{Z_3(J)} }.
\end{multline*}
As the operator $B_2$ is of similar type, we are able to deal with the difference
\begin{align*}
B_2(w_1,\rho_1) c_1 - B_2(w_2,\rho_2) c_2 = [B_2(w_1,\rho_1) - B_2(w_0,\rho_0)] c
+ [B_2(w_1,\rho_1) - B_2(w_2,\rho_2)] c_2
\end{align*}
in the same way.

We finally consider the difference $B_3(w_1,\rho_1)\rho_1 - B_3(w_2,\rho_2)\rho_2$.
Using again \eqref{diff:nl} and the embeddings, we may proceed as follows
\begin{align*}
\norm{ B_3(w_1,\rho_1)\rho_1 - B_3(w_2,\rho_2)\rho_2}{ X_1(J) }
\le
\norm{ [B_3(w_1,\rho_1) - B_3(w_2,\rho_2)] \rho_1}{ X_1(J) }
\\
{}\quad+\norm{ B_3(w_2,\rho_2) \rho }{ X_1(J) }
\le
C_1 \norm{ (\rho,c) }{ Z_\mu(J) }
+
C_2 \norm{ \rho }{ Z_\mu(J) }
\\
\le
C_1 T^{1/2}  \norm{c}{ \zero{Z_2(J)} } + (C_1 + C_2) ( T^{3/4 + 1} + T^{1/2 + 1}) \norm{\rho}{ \zero{Z_3(J)} }.
\end{align*}
These estimates show that there exists a $\kappa_3(T)$ independent of $w'$ and $w_i$, $i=1,2$, 
and with the property $\kappa_3(T) \to 0$ as $T \to 0$ such that
\begin{equation*}
\norm{(u',c',\mu')}{ \zero{Z_1(J)} \times \zero{Z_2(J)} \times \zero{Z_\mu(J)} }
\le \kappa_3(T) \norm{ (u,c,\mu) }{  \zero{Z_1(J)} \times \zero{Z_2(J)} \times \zero{Z_\mu(J)} }.
\end{equation*}
{\bf Step 2: Self-mapping.} We may proceed similarly as in the previous part, however we are now concerned with
higher regularity estimates. To begin with,
let $w=(u,c,\mu) \in \Sigma_T$ be given. We have to show that $w'=(u',c',\mu')$, given as the solution of
$\mathcal{L}w' = \mathcal{F}(w,\rho)$ with $\rho = L[u]\rho_0$, lies in $\Sigma_T$ as well,
that is $\norm{ w' - \ol{w} }{ \zero{\mathcal{Z}_{1,\mathcal{B}}(J) \times
 \zero{\mathcal{Z}_2(J)} } \times \zero{\mathcal{Z}_\mu(J)} } \le 1$. By the Theorems \ref{theo:ch} and \ref{theo:stokes:3}
the following estimate is available.
\begin{multline*}
\norm{ w' - \ol{w} }{\zero{\mathcal{Z}_{1,\mathcal{B}}(J)} \times \zero{ \mathcal{Z}_2}(J) \times \zero{\mathcal{Z}_\mu(J)}  }
\le
M \left \{
\norm{ F_1(w,\rho) - F_1(w_0,\rho_0)}{ \zero{\mathcal{X}_{1,n,\Gamma}}(J) } \right.
\\
+\norm{F_2(w,\rho) -F_2(w_0,\rho_0) }{ \mathcal{X}_2(J) } + \norm{ F_\mu(w,\rho)-F_\mu(w_0,\rho_0)}{
\zero{\mathcal{X}_\mu(J)} }
\left. \right \}
\\
\le M \big \{
\norm{ B_{low}(w,\rho) - B_{low}(w_0,\rho_0)}{ \zero{\mathcal{X}_{1,n,\Gamma}}(J) }
+ \norm{ B_{1}(w,\rho)u }{ \zero{\mathcal{X}_{1,n,\Gamma}}(J) }\\
+\norm{ B_{2}(w,\rho)c }{ \zero{\mathcal{X}_{1,n,\Gamma}}(J) }
+\norm{ B_{\mu}(w,\rho)\mu }{ \zero{\mathcal{X}_{1,n,\Gamma}}(J) }
+\norm{ B_{3}(w,\rho)\rho- B_{3}(w_0,\rho_0)\rho_0}{ \zero{\mathcal{X}_{1,n,\Gamma}}(J) }\\
+
\norm{F_2(w,\rho) -F_2(w_0,\rho_0) }{ \mathcal{X}_2(J) } + \norm{ F_\mu(w,\rho)-F_\mu(w_0,\rho_0)}{
\zero{\mathcal{X}_\mu(J)} }
\big \}.
\end{multline*}
We begin with estimating the differences of $F_2$ and $F_\mu$.
By means of the estimate from Lemma \ref{lem:prod:2}, combined with \eqref{diff:nl} and
\eqref{eq:sob-T} we obtain
\begin{multline*}
\norm{F_2(w,\rho) - F_2(w_0,\rho_0) }{ \mathcal{X}_2(J) }
\le
\norm{ \tfrac{\varepsilon_0}{\gamma_0} }{ \cont{}{}{\ol{\Omega}} }
\big(
\norm{ \rho_0 - \rho }{  \cont{}{}{J;\cont{}{}{\ol{\Omega}}} }
\norm{ \del{}{t} c }{ \mathcal{X}_2(J) }
\\
+\norm{ \del{}{t} \rho }{  \mathcal{X}_2(J) } \norm{ c }{ \cont{}{}{J;\cont{}{}{\ol{\Omega}}} }
+
T^{1/p} \norm{ \rho u c - \rho_0 u_0 c_0 }{ \cont{}{}{J;\cont{1}{}{\ol{\Omega};\R^n} } }\\
+
\norm{ \gamma_0 - \gamma}{ \cont{}{}{J;\cont{}{}{\ol{\Omega}}} }
\norm{ \mu }{ \mathcal{Z}_\mu(J) }
+
\norm{ \gamma_0 - \gamma }{ \lp{p}{J;\cont{1}{}{ \ol{\Omega} }} }
\norm{ \mu }{ \cont{}{}{J;\sob{1}{p}{\Omega}} }
\big)\\
+
\norm{ \tfrac{\varepsilon^2_0}{\gamma_0} }{ \cont{}{}{\ol{\Omega}} }
\norm{ \tfrac{\varepsilon_0}{\gamma_0} }{ \cont{1}{}{\ol{\Omega}} }
\norm{ \nabla [\mu - \mu_0] }{ \mathcal{X}_2(J)}
\\
\le  C_1 T + C_2 T^{1/p}  + C_3 T^{1/p} + C_4 [T^{1/2+\beta} + T] + C_5[ T^{1/p+1/4+\beta} + T^{1/p + 1/2}] + C_6 T^{1/4}
\end{multline*}
and
\begin{multline*}
\norm{F_\mu(w,\rho) - F_\mu(w_0,\rho_0) }{ \mathcal{X}_\mu(J) }
\le C_1 (
\norm{ \varepsilon_0 - \varepsilon}{ \zero{\mathcal{Z}_2(J)} }
\norm{\Delta c}{ \cont{}{}{J \times \ol{\Omega} } }
\\
+\norm{ \varepsilon_0 - \varepsilon}{ \cont{}{}{J \times \ol{\Omega} } }
\norm{\Delta c}{ \mathcal{Z}_2(J) }
+
\norm{ \varepsilon_0 - \varepsilon }{ \lp{p}{J;\cont{1}{}{\ol{\Omega}}} }
\norm{\Delta c}{ \cont{}{}{J;\sob{1}{p}{\Omega}}}
)
\\
+C_2\norm{ (\rho,\nabla \rho, c, \nabla c)  - (\rho_0, \nabla \rho_0, c_0, \nabla c_0) }{ \zero{\mathcal{Z}_2(J)} }\\
\le
C_1 \big \{   C_{11} [ T^{1/2} + T^{1/p+1}] + C_{12} [ T^{1/2+\beta} + T]
+
C_{13}[ T^{1/p + 1/4 + \beta} + T^{1/p + 1/2}] \big \}\\
+
C_2 \big \{ T^{1/2} + T^{1/p + 1} + T^{1/4} + T^{1/2+1/p} + T^{1/p} \big \}.
\end{multline*}
We give some details for the treatment of one of the critical terms involving a vector-valued
Bessel potential space, where we have to use Lemma \ref{lem:prod:2}. 
In order to estimate $\norm{F_\mu(w,\rho) - F_\mu(w_0,\rho_0) }{ \mathcal{X}_\mu(J) }$ 
we have to estimate, among others, the term ${(\varepsilon - \varepsilon_0)
\Delta c}$ in the space ${}_0\sob{1/2}{p}{J;\lp{p}{\Omega}}$. 
Setting $m=\varepsilon - \varepsilon_0$ we proceed as follows.
\begin{align*}
\norm{m
\Delta c}{{}_0\sob{1/2}{p}{J;\lp{p}{\Omega}}} &\le 
\norm{m
(\Delta c-\Delta \bar{c})}{{}_0\sob{1/2}{p}{J;\lp{p}{\Omega}}}+
\norm{m
(\Delta \bar{c}-\Delta c_0)}{{}_0\sob{1/2}{p}{J;\lp{p}{\Omega}}}
\\
& \quad +\norm{m
\Delta c_0}{{}_0\sob{1/2}{p}{J;\lp{p}{\Omega}}}
\\
& \le C\Big(\norm{m}{\cont{}{}{J\times\bar{\Omega}}}\big[\norm{\Delta c-\Delta \bar{c}}{{}_0\sob{1/2}{p}{J;\lp{p}{\Omega}}}
+
\norm{\Delta \bar{c}-\Delta {c}_0}{{}_0\sob{1/2}{p}{J;\lp{p}{\Omega}}}\big]
\\
& \quad +\norm{m}{\cont{1/2+\epsilon}{}{J;\cont{}{}{\bar{\Omega}}}}
\big[\norm{\Delta c-\Delta \bar{c}}{\lp{p}{J\times \Omega}}+\norm{\Delta \bar{c}-\Delta {c}_0}{\lp{p}{J\times \Omega}}\big]
\\
& \quad +\norm{m \Delta c_0}{\cont{1/2+\epsilon}{}{J;\cont{}{}{\bar{\Omega}}}}\Big);
\end{align*}
here $\epsilon>0$ can be chosen sufficiently small so that $m$ belongs to
the required spaces. Using the definition of $\Sigma_T$, the structure of $m$, and
additional regularity of $m$, it is then not difficult to see that each of the five
terms can be estimated by a quantity that is independent of $c$ (and $\rho$) and tends to $0$
as $T\to 0$.

To obtain a similar estimate for differences of $B_{low}$ and the $B_{i}$-terms, $i=1,2,\mu,3$, we are able
to proceed in the same way as above, since $\mathcal{X}_{1,n,\Gamma}(J)$ coincides with $\mathcal{X}_\mu(J)^n$
for functions with vanishing initial traces, and the terms occurring in these differences
have the same structure as above, namely higher order terms multiplied by a 'small' difference
as well as lower order terms with more time regularity. Therefore, the desired estimates
for the remaining terms can be derived similarly as above, so that we will omit the details here.

Finally, putting together all estimates from above we obtain an inequality of the form
\begin{align*}
\norm{ (u',c',\mu') - (\ol{u},\ol{c},\ol{\mu})}{ \zero{\mathcal{Z}_1(J)} \times
\zero{\mathcal{Z}_2(J)} \times \zero{\mathcal{Z}_\mu(J)} } \le \kappa(T)
\end{align*}
with $\kappa(T)$ independent of $w, w'$ and tending to $0$ as $T \to 0$. Therefore, choosing $T$ so small
that $\kappa(T) \le 1$, the fixed point mapping enjoys the self-mapping property.

To conclude, we have seen that for sufficiently small $T$ the mapping
$\mathcal{G}: \Sigma_T \mapsto \Sigma_T$ is a strict contraction w.r.t.\
the topology of $Z(J)$, hence by Lemma \ref{lem:1} and the contraction mapping principle
admits a unique fixed point $(u,c,\mu)$ in $\mathcal{Z}_1(J) \times \mathcal{Z}_2(J)
\times \mathcal{Z}_\mu(J)$. This fixed point together with $\rho = L[u] \rho_0 \in \mathcal{Z}_3(J)$
yields a unique local in time strong solution $(u,c,\rho)\in \mathcal{Z}(J)$ of problem \eqref{eq:nsch-1}-\eqref{bc:3}.
This solution can be extended by the standard method of successively repeating
the above arguments on
intervals $[t_i , t_{i+1} ]$. Either after finitely many steps we reach $T_0$,
or we have an infinite strictly increasing sequence which converges to some
$T^*(u_0,c_0,\rho_0) \le T_0$. In case $\lim_{i \to \infty} (u,c,\rho)(t_i)
=: (u(T^*),c(T^*),\rho(T^*))$ exists in the phase space $\mathcal{V}_p$,
we may continue the process, which shows that the maximal time is characterized by
condition \eqref{cond:blow-up}.
This completes the proof of Theorem \ref{theo:main:1}.
\eproof

\section{The nonlinear problem with inhomogeneous boundary conditions} \label{sec:8}
In this section we will discuss the case of inhomogeneous boundary conditions and how these lead to some restrictions.
We now consider \eqref{eq:nsch-1}-\eqref{ic} with the following inhomogeneous boundary
conditions, that is we replace \eqref{bc:1}-\eqref{bc:3} by
\begin{equation} \label{bc:4}
\begin{split}
u & = h_d(t,x), \quad (t,x) \in J \times \Gamma_d, \quad
\\
\ska{ u }{ \nu }{|\Gamma_s} & = h_{s1}(t,x), \quad
\mathcal{Q}  \mathcal{S} \cdot \nu \equiv 2 \eta \mathcal{Q} \mathcal{D}(u) \cdot \nu_{|\Gamma_s} = h_{s2}(t,x),
\quad (t,x) \in J \times \Gamma_s,
\\
\del{}{\nu} c & = h_1(t,x), \quad \del{}{\nu} \mu = h_\mu(t,x),
\quad (t,x) \in J \times \Gamma,
\end{split}
\end{equation}
where $h_d$ and $h_s := (h_{s1},h_{s2})$ have to be subject to the conditions
\begin{align*}
\ska{ h_d(t,x) }{ \nu(x) }{|\Gamma_d} \ge 0, \quad \forall (t,x) \in J \times \Gamma_d,
\quad
h_{s1}(t,x) \ge 0, \quad \forall (t,x) \in J \times \Gamma_s,
\end{align*}
in order to ensure that Lemma \ref{lem:rho} is applicable to solve the continuity equation \eqref{eq:me}.

The homogeneous boundary conditions were very convenient to find the estimate
\eqref{est:diff-ch}, that is to estimate $(c',\mu')$ in a weaker topology. In fact, to get \eqref{est:Phi-1}
we used the fact that certain boundary integrals vanish in view of the homogeneous boundary
conditions. This is no longer true, however, those additional terms for which we cannot argue as above
are of lower order, so that we are able to estimate them 'directly' instead of integrating by parts
to get rid of one spatial derivative. For instance, one has to deal with
the term $\du{ \tfrac{\varepsilon_0}{\gamma_0} A_{\gamma_1 - \gamma_2} \mu_2}{ \psi }$ where
$\norm{H^*\psi}{ Y } \le 1$. H\"older's inequality yields
\begin{align*}
|\du{ \tfrac{\varepsilon_0}{\gamma_0} A_{\gamma_1 - \gamma_2} \mu_2 }{ \psi }|  \le &\,
C(\varepsilon_0,\gamma_0) \norm{\psi}{ \lp{2}{J;\lp{2}{\Omega}} } \norm{ A_{\gamma_1 - \gamma_2} \mu_2 }{ \lp{2}{J;\lp{2}{\Omega}} }
\\
\le &\,  C T^{1/4} \big(
\norm{ \gamma_1 - \gamma_2 }{ \lp{p_1}{J;\lp{p_1}{\Omega}} } \norm{ \Delta \mu_2}{ \lp{p}{J;\lp{p}{\Omega}} }
\\
&\, +
\norm{ \nabla[ \gamma_1 - \gamma_2] }{ \lp{p_1}{J;\lp{2}{\Omega;\R^n}} } \cdot
\norm{ \nabla \mu_2}{ \lp{p}{J;\cont{1}{}{\ol{\Omega};\R^n}} }
\big)\\
\le & \,C T^{1/4} \big( C_1 \norm{ (c,\rho) }{ \lp{p_1}{J;\lp{p_1}{\Omega}} }
+
C_2 \norm{ \nabla (\rho,c) }{ \lp{p_1}{J;\lp{2}{\Omega;\R^n}} } \big)
\\
\le & \,C T^{1/4} \big ( \norm{c}{ \zero{Z_2(J)} } + \norm{ \rho }{ \zero{Z_3(J)} } ),
\end{align*}
since $\sob{1}{2}{\Omega} \hookrightarrow \lp{p_1}{\Omega}$, $p_1 = 2p/(p-2)$.

Another problem arises from the boundary condition
$2 \eta \mathcal{Q} \mathcal{D}(u) \cdot \nu = h_{s2}$ on $\Gamma_s$, which becomes
nonlinear when $\eta = \eta(\rho,c)$. According to Theorem \ref{theo:stokes:3} this boundary
equation has to be considered in $\mathcal{Y}_{1,s}(J;\R^n) : = \sobb{\frac{3}{2}-\frac{1}{2p}}{p}{J;\lp{p}{\Gamma_s;\R^n}} \cap
\lp{p}{J;\sobb{3-\frac{1}{p}}{p}{\Gamma_s;\R^n}}$, in particular, the coefficient $\eta(\rho,c)$ has to be
in this space as well (this space forms a multiplication algebra for $p > \hat{p}$). However, looking for
$c \in \mathcal{Z}_2(J)$ and $\rho \in \mathcal{Z}_3(J)$, we obtain by trace theory that
\begin{align*}
& c_{|\Gamma_s} \in \sobb{1-\frac{1}{4p}}{p}{J;\lp{p}{\Gamma_s}} \cap \lp{p}{J;\sobb{4-\frac{1}{p}}{p}{\Gamma_s}},
\\
& \rho_{|\Gamma_s} \in \sobb{2+\frac{1}{4} -\frac{1}{3p}}{p}{J;\lp{p}{\Gamma_s}} \cap \cont{}{}{J;\sobb{3-\frac{1}{p}}{p}{\Gamma_s}},
\end{align*}
and comparing these regularities with $\mathcal{Y}_{1,s}(J;\R^n)$, one perceives that $c_{|\Gamma_s}$ does not possess enough time
regularity. Hence, permitting non-zero boundary data $h_{s2}$ would
lead to the restriction $\eta=\eta(\rho)$. But from the physical point of view it is more reasonable to consider
viscosities depending on $c$ to model different viscosities for different phases.
We therefore abstain from considering boundary data $h_{s2} \not = 0$ in the theorem below.

The subsequent theorem is an extension of Theorem \ref{theo:main:1} that takes into account inhomogeneous boundary data.
\begin{theorem} \label{theo:main:2}
Let $\Omega$ be a bounded domain in $\R^n$, $n \ge 1$, with
compact $C^{4}$-boundary $\Gamma$ decomposing disjointly as
$\Gamma = \Gamma_d \cup \Gamma_s$ with dist$\,(\Gamma_d,\Gamma_s)>0$, $J_0=[0,T_0]$ with $T_0 \in
(0,\infty)$, and $p \in (\hat{p},\infty)$. Assume further that the
following conditions are satisfied.
\begin{compactenum}
\item[(i)] $\varepsilon\in \cont{4}{}{\R^2}$, $\ol{\psi}\in \cont{5}{}{\R^2}$, $\gamma\in \cont{2}{}{\R^2}$;
$\eta,\,\lambda\in \cont{4}{}{\R^2}$;
\item[(ii)] $\eta$, $2\eta+\lambda$, $\varepsilon$, $\gamma>0$ in $\R^2$;
\end{compactenum}
Then for each $f_{ext} \in \mathcal{X}_1^n(J_0)$
and initial data $w_0=(u_0,c_0,\rho_0)$ in
\begin{align*}
\mathcal{V} := \sobb{4 - \frac{2}{p} }{p}{\Omega;\R^n} \times
\sobb{4 - \frac{4}{p} }{p}{\Omega} \times
\{ \varphi \in \sob{3}{p}{\Omega;\R_+}: \varphi(x) > 0, \quad \forall x \in \ol{\Omega} \}
\end{align*}
and boundary data
\begin{align*}
h_d & \in \mathcal{Y}_{0,d}(J;\R^n), \quad \ska{h_d}{\nu}{} \ge 0, \quad
h_s = (h_{s1},h_{s2}) \in \mathcal{Y}_{0,s}(J) \times \{0\}, \quad
h_{s1} \ge 0,
\\
h_1 & \in \mathcal{Y}_1(J), \quad h_\mu \in \mathcal{Y}_\mu(J),
\end{align*}
satisfying the compatibility conditions (CP)
\begin{equation*}
\begin{split}
& u_{0|\Gamma_d} = h_{d|t=0} \in \sobb{4-3/p}{p}{\Gamma_d;\R^n},
\quad
\ska{ u_0 }{ \nu }{|\Gamma_s} = h_{s1|t=0} \in \sobb{4-3/p}{p}{\Gamma_s},
\quad
\mathcal{Q} \mathcal{S}_{|t=0} \cdot \nu_{|\Gamma_s} = 0,
\\
& \del{}{\nu} c_{0|\Gamma} = h_{1|t=0} \in \sobb{3-5/p}{p}{\Gamma},
\quad
\del{}{\nu} \mu(\rho_0,c_0)_{|\Gamma} = h_{2|t=0} \in \sobb{1-5/p}{p}{\Gamma},
\\
& \rho_{0|\Gamma_d} \del{}{t} h_{d|t=0} - \divv \mathcal{S}_{|t=0,\Gamma_d}
= ( \divv \mathcal{P}_{|t=0} - \rho_0 \nabla u_0 \cdot u_0
+ \rho_0 f_{ext|t=0} )_{|\Gamma_d} \in \sobb{2-\frac{3}{p}}{p}{\Gamma_d;\R^n},
\\
& \rho_{0|\Gamma_d} \del{}{t} h_{s1|t=0}
- \ska{ \divv \mathcal{S}_{|t=0} }{ \nu }{|\Gamma_s}
= \ska{ \divv \mathcal{P}_{|t=0} - \rho_0 \nabla u_0 \cdot u_0
+ \rho_{0} f_{ext|t=0} }{ \nu }{|\Gamma_s}
\in \sobb{2 - \frac{3}{p}}{p}{\Gamma_s},
\\
& \del{}{t} h_{s2|t=0} - \mathcal{Q} \mathcal{S}( \rho_0^{-1} \divv \mathcal{S} )_{|t=0,\Gamma_s}
\cdot \nu_{|\Gamma_s}
\\
& \hspace{2cm}
=\mathcal{Q} \mathcal{S}( \rho_0^{-1} \divv \mathcal{P}_{|t=0}
-  \nabla u_0 \cdot u_0 + f_{ext|t=0})_{|\Gamma_s}
\cdot \nu_{|\Gamma_s} \in \sobb{1 - \frac{3}{p}}{p}{\Gamma_s;\R^n},
\end{split}
\end{equation*}
there is a unique solution $w=(u,c,\rho)$ of \eqref{eq:nsch-1}-\eqref{ic}, \eqref{bc:4}
on a maximal time interval $J_*:=[0,T^*)$, $T^* \le T_0$, if the solution is not global;
the solution $w$ belongs to the class $\mathcal{Z}(J_1)$ for each interval
$J_1 = [0,T_1]$ with $T_1 < T^*$, or to the class $\mathcal{Z}(J_0)$ if the solution exists globally.
The maximal time interval $J_*$ is characterized by the property:
\begin{align*}
\lim_{t \to T^*} w(t) \quad \text{does not exist in $\mathcal{V}_p$,}
\end{align*}
where $\mathcal{V}_p$ is defined as the space of all $(u_1,c_1,\rho_1)\in \mathcal{V}$ such that the
compatibility conditions (CP) hold with $(u_0,c_0,\rho_0)$ being
replaced by $(u_1,c_1,\rho_1)$.
Moreover, for fixed $f_{ext}$, $h_d$, $h_s$, $h_1$, $h_\mu$ not depending on $t$
the solution map $w_0 \rightarrow w(\cdot)$ generates a local semiflow on
the phase space $\mathcal{V}_{p}$.
\end{theorem}
\end{document}